\documentclass[11pt]{amsart}

\allowdisplaybreaks[2] \textwidth15.1cm \textheight22cm \headheight12pt \oddsidemargin.4cm
\input{xy}
\xyoption{all}
\usepackage{graphicx}
\usepackage{color}
\usepackage{verbatim}
\usepackage[bookmarks,colorlinks,breaklinks]{hyperref}  
\hypersetup{linkcolor=blue,citecolor=red,filecolor=blue,urlcolor=red} 
\usepackage{amssymb,amsmath,amsthm,amsfonts}
\usepackage{enumitem}
\usepackage{setspace,bbm}
\usepackage{empheq}
\usepackage[all]{xy}
\usepackage[normalem]{ulem}  
\usepackage{chngcntr}

\newtheorem{theorem}{Theorem} [section]
\newtheorem{prop}[theorem]{Proposition}
\newtheorem{lemma}[theorem]{Lemma}
\newtheorem{cor}[theorem]{Corollary}
\newtheorem{conjecture}[theorem]{Conjecture}
\newtheorem{proposition}[theorem]{Proposition}
\newtheorem{definition}[theorem]{Definition}
\newtheorem{question}[theorem]{Question}

\theoremstyle{remark}
\newtheorem{remark}{Remark}

\numberwithin{equation}{section}
\numberwithin{figure}{section}

\newcommand\A{{\mathbb A}}
\newcommand\C{{\mathbb C}}
\newcommand\CC{{\mathbb C}}

\newcommand\N{{\mathbb N}}

\renewcommand\P{{\mathbb P}}
\newcommand\PP{{\mathbb P}}
\newcommand\R{{\mathbb R}}
\newcommand\RR{{\mathbb R}}
\newcommand\Z{{\mathbb Z}}

\newcommand\Q{{\mathbb Q}}

\renewcommand\phi{\varphi}

\newcommand\cM{\mathcal{M}}

\newcommand\Gal{\operatorname{Gal}}

\newcommand\supp{\operatorname{supp}}   

\renewcommand\gcd {\operatorname{gcd}} 
\newcommand\capacity {\operatorname{Cap}} 
\newcommand\Res {\operatorname{Res}} 

\newcommand\Per {\mathrm{Per}}

\newcommand\Sing {\mathrm{Sing}}

\newcommand\al {\alpha}
\newcommand\la {\lambda}

\newcommand\kbar {\overline{k}}
\newcommand\Qbar {\overline{\Q}}
\newcommand\kvbar {\overline{k}_v}
\newcommand\ksep {k^{\rm sep}}

\newcommand\ord {\operatorname{ord}}
\newcommand\hhat {\hat{h}}

\definecolor{myblue}{rgb}{0.6, 0.9, 1}
\definecolor{mygreen}{rgb}{0,0,1}
\definecolor{purple}{rgb}{0.6,0.2,1}
\definecolor{orange}{rgb}{0.8,0,0.2}

\begin{document}

\title{Quasi-adelic measures and equidistribution on $\P^1$}

\author{Niki Myrto Mavraki}
\thanks{Email: myrtomav@math.ubc.ca, Dept. of Mathematics, University of British Columbia, Vancouver, BC V6T 1Z2, Canada, Telephone: +1-604-728-7631}

\author{Hexi Ye}
\thanks{Email: yehexi@gmail.com, Dept. of Mathematics, Zhejiang University, Hangzhou, 310027, P.R.China, Telephone: +86-571-8795-2335}


\begin{abstract}
Baker--Rumely \cite{Baker:Rumely:equidistribution,BRbook}, Favre--Rivera-Letelier \cite{FRL:equidistribution} independently proved an important arithmetic equidistribution theorem for points of small height on the Berkovich compactification of the projective line with respect to an adelic measure on $\P^1$. We generalize the notion of an adelic measure to that of a quasi-adelic measure on $\PP^1$, and show that arithmetic equidistribution of points with small height holds for quasi-adelic measures as well. Moreover, we show that the canonical measure associated with a dynamical pair $(f,c)$ on $\P^1$ is rarely adelic. 
We prove that for certain examples of families of rational functions parameterized by $\PP^1$, corresponding to the curve $\Per_1(\la)$ introduced by Milnor for a root of unity $\la$, the measure corresponding to a general starting point is quasi-adelic. Finally, we place our results in context by establishing their connection with two problems in arithmetic dynamics.
\end{abstract}

\maketitle
\tableofcontents
\thispagestyle{empty}

\section{Introduction}
Baker--Rumely \cite{Baker:Rumely:equidistribution, BRbook} and  Favre--Rivera-Letelier \cite{FRL:equidistribution}, relying on a potential-theoretic approach, independently proved an important arithmetic equidistribution theorem for points of small height on the Berkovich compactification of the projective line. 
Around the same time, Chambert-Loir \cite{CLoir}, using an approach based on Arakelov theory, proved a more general version of this arithmetic equidistribution theorem in the setting of curves. 
Let us now describe the arithmetic equidistribution theorem. 
 Let $k$ be a product formula field, for example a number field or the function field of a smooth and projective curve, and let $\mu=\{\mu_v\}_{v\in\cM_k}$ be a collection of probability measures $\mu_v$ on the Berkovich projective line $\P^{1, an}_{v}$, one for each place $v$ of $k$. 
We say that $\mu$ is an adelic measure on $\P^1$ if all the measures $\mu_v$ have continuous potentials \textbf{and} those potentials are trivial at all but finitely many places of $k$. 
To each adelic measure we can associate a height function given by the sum of potential functions of the measures ${\mu_v}$. 
For this choice of height function, the main theorem in \cite{BRbook, FRL:equidistribution} asserts that points of small height equidistribute with respect to $\mu_v$ at all places $v$ of $k$. We remark that the `adelic' hypothesis appears in the equidistribution results on more general varieties  other than $\P^1$ \cite{CLoir, Thuillier:these, Yuan}.  
There it is subsumed in the notion of an adelic metrized line bundle \cite{Zhang:line, Zhang:metrics}. 

In this article, we introduce the notion of a quasi-adelic measure $\mu=\{\mu_v\}$ on $\P^1$, which extends the notion of an adelic measure as we do not require the potentials to be trivial at all but finitely many places. Instead, we  impose a certain summability condition. As in the case of an adelic measure, there is a natural height function defined on $\PP^1(\ksep)$, which we denote by $\hat{h}_{\mu}$, associated to each quasi-adelic measure; see \S\ref{Quasi-adelic measure and canonical height function}. 
In contrast to the height associated with an adelic measure, this height allows non-trivial contributions from infinitely many local heights. 
We also point out that we can uniquely define a metric on $\mathcal{O}_{\P^1}(1)$ associated with a quasi-adelic measure. 
In this quasi-adelic metrized line bundle, our summability condition controls the total distortion of this metric from the trivial metric over all places $v$ of $k$, thereby making it more flexible than the adelic metric in \cite{Zhang:line, Zhang:metrics} as there can be infinitely places $v\in\cM_k$ at which the metric $\|\cdot\|_v$ is non-trivial. Our first main result extends that in \cite{BRbook} and \cite{FRL:equidistribution} to the setting of quasi-adelic measures, building upon the potential-theoretic approach therein. 
\begin{theorem}
\label{arithmetic equidistribution}
Let $k$ be a product formula field and ${\mathbb \mu} = \{ \mu_v \}_{v \in M_k}$ be a quasi-adelic measure.  Suppose that $S_n$ is a sequence of $\Gal(\ksep/k)$-invariant subsets of $\PP^1(\ksep)$ such that $|S_n|\to \infty$ and $\hhat_{{\mathbb \mu}}(S_n) \to 0$
as $n \to \infty$. Then for each $v \in \cM_k$,
 the sequence of probability measures $[S_n]_v$, weighted equally on the points in $S_n$,  converges weakly to $\mu_v$ on $\PP^{1, an}_{,v}$ as $n \to \infty$. 
 \end{theorem}

Our next main theorem shows that many measures which arise naturally in applications to the study of dynamical systems are not adelic. The measures we study come from dynamical pairs $(f,c)$ over $\P^1$, where $f\in k(t)(z)$ with $\deg f\ge 2$ and $c\in k(t)$ is a starting point.
To each non-isotrivial and non-preperiodic dynamical pair we associate a \emph{canonical measure} $\mu_{f,c}=\{\mu_{f,c,v}\}_{v\in\cM_k}$; see \S\ref{dynamicalmeasure}. Let us point out that when $f\in k(z)\subset k(t)(z)$ has degree $d\geq 2$ and $c(t)=t$, this measure is the one studied in \cite{BH:polynomials, Baker:Rumely:equidistribution}. 
 If $\mu_{f,c}$ is adelic (respectively quasi-adelic), we call the dynamical pair $(f,c)$ adelic (respectively quasi-adelic). 
We postpone the statement of our theorem until Section \S\ref{a pair is rarely adelic} (see Theorem~\ref{genericnoadelic}) when all the notation will be in place. We state here a corollary instead. 

\begin{theorem}\label{cornoadelic}
Let $k$ be a number field or the function field of a smooth projective curve defined over a field of characteristic zero, and consider $f\in k(t)(z)$ with $d:=\deg f\ge 2$. Assume that there is an $N\in\N$ and an $\al\in\P^1(\overline{k})$ such that for $g:=f^N$, we have $2\le\deg(g_{\al})<d^N$ and $g^2_{\al}$ is not a polynomial map. Then there is a constant $L\ge0$ such that for all $c\in k(t)$ with $\hat{h}_{g_{\al}}(c(\al))>L$, the dynamical pairs $(f,c)$ are not adelic.
\end{theorem}
The conditions in our theorem are in fact mild and are satisfied ``generically''. If $k'$ is a number field such that $c\in k'(t)$ and $\al\in k'$, then the Northcott property of heights gives that inequality $\hat{h}_{g_{\al}}(c(\al))>L$ will be satisfied for all but finitely many $c(\al)\in k'$. Our assumption that $\deg(g_{\al})<d^N$ means that $t=\al$ is a point of bad reduction for $g_t$. 
That is, if $g_t=P_t/Q_t$ with $P_t, Q_t\in k(t)(z)$ of degree $d^N\ge 3$ having no common factor in $\overline{k(t)}$, then we have that the specialized map $g_{\al}$ has degree at most $d^N-1$ if and only if the resultant of $P_t$ and $Q_t$ vanishes at $t=\al$. Moreover, the map $g_{\al}$ has degree at most one, if and only if $P_{\al}$ and $Q_{\al}$ share at least $d^N-1$ common factors in $\overline{k}$.

It is now natural to ask whether the dynamical measures that are not adelic are quasi-adelic. The first example of quasi-adelic measures was studied by DeMarco, Wang and Ye  \cite{DWY:Per1}. 
They considered the family of rational maps $g_{\lambda,t}(z)=\frac{\lambda z}{z^2+tz+1}$ for $\lambda\in \overline{\Q}\backslash \{0\}$, which corresponds to the family of conjugacy classes  $\Per_1(\lambda)$ in the moduli space of quadratic maps introduced by Milnor \cite{Milnor:quad}. They showed that if the starting point is $c\in\{1,-1\}$, which is a critical point of $g_{\lambda,t}$, and if $\lambda\neq 0$ is not a root of unity of order at least $2$, then the dynamical pair $(g_{\lambda}, c)$ is quasi-adelic, but not adelic. 
Using the equidistribution Theorem  \ref{arithmetic equidistribution} proved in this article,
 they proved that when $\la\neq0$ the curve $\Per_1(\lambda)$ contains at most finitely many postcritically-finite maps. 
Note that if $\lambda$ is not a root of unity of order $\ell\geq 2$, then any $n$-th iterate of $g_{\lambda, t}$ degenerates to the same constant map at $t=\infty$. 
However, if $\lambda$ is a root of unity of order $\ell\geq 2$, although the $n$-th iterate $g^n_{\lambda, t}$ degenerates to a constant map when $n< \ell$, the $\ell-$th iterate $f_t:=g^{\ell}_{\lambda, t}$ degenerates to the degree $2$ map $f_\infty(z)=z/(z^2+1)$, which is not the $\ell$-th iterate of $g_{\lambda, \infty}$; see Proposition \ref{degenerations}. We refer the reader to \cite{Eps,D:moduli} for more about this phenomenon. 
In this article, we extend the methods from \cite{DWY:Per1} and study families of rational maps corresponding to $\Per_1(\la)$ for a root of unity $\la$. In doing so, we provide more examples of quasi-adelic measures that fail to be adelic.
\begin{theorem}\label{mainthm3}
Let $\la$ be a primitive  $\ell$-th root of unity with $\ell\geq 2$ and $k$ be a number field with $\la\in k$. Consider the rational function $g_{\la,t}(z)=\frac{\la z}{z^2+tz+1}\in k(t)(z)$ and let $c(t)\in k(t)$ be such that $0$ and $\infty$ are not in the orbit of $c(\infty)$ iterated under the map $f_{\infty}(z)=z/(z^2+1)$.
Then the measure $\mu_{g_{\la},c}=\{\mu_{g_{\la},c,v}\}_{v\in\cM_k}$ is quasi-adelic. Furthermore, if  $c(\infty)$ is not a preperiodic point of $f_{\infty}$, then the measure $\mu_{g_{\la},c}$ is not adelic. 
\end{theorem}
We emphasize here that our assumptions on $c(\infty)$, namely that $0$ and $\infty$ are not in the orbit of $c(\infty)$ under the iteration by $f_{\infty}$ and that $c(\infty)$ is not a preperiodic point for $f_{\infty}$, are very mild.  Indeed, all we are assuming is that $c(\infty)\notin\mathcal{P}$ for a set $\mathcal{P}$ of bounded Weil height. 
In particular, from the Northcott property of heights, we have that for any $\kappa\in\N$ there are only finitely many choices for $c(\infty)$ in $\mathcal{P}$ with $[\Q(c(\infty)):\Q]\le \kappa$. As a special case, our Theorem \ref{mainthm3} holds when the starting point $c\in\{1,-1\}$ is a critical point of $g_{\la,t}$.

Our next theorem explains the connection between our definition of a quasi-adelic measure and two classical problems in arithmetic dynamics that motivated our research. 
We first state our result and subsequently provide context. 
Let  $K=k(X)$ be the function field of a smooth projective curve $X$ defined over a number field $k$, and consider $f\in K(z)$ with degree $d\ge 2$ and $c\in K$. 
We have three heights associated with the dynamical pair $(f,c)\in K(z)\times K$; the Weil height $h:X(\overline{k})\to \R_{\ge 0}$ associated to a degree $1$ divisor on $X$, the Call-Silverman canonical height $\hat{h}_f(c)$ of $c$ associated with the map $f$ defined over the function field $K$, and for each $t\in X(\overline{k})$ such that the specialized map $f_t$ has degree at least $2$, the Call-Silverman canonical height $\hat{h}_{f_t}(c(t))$. 
We consider the case $X=\P^1$ and may take $h(t)$ to be the Weil height associated with the divisor $([1:0])\in\mathrm{Pic}(\P^1)$.
\begin{theorem}\label{vchunlikely}
\counterwithin{enumi}{theorem}
Let $k$ be a number field. Let $f\in k(t)(z)$ and $c\in k(t)$ be such that the degree $\deg_zf\ge 2$ and the dynamical pair $(f,c)$ is quasi-adelic. The following hold.
\begin{enumerate}
\item  \label{thm3part2}
As $t\in\P^1( \overline{k})$ varies we have
$\hat{h}_{f_{t}}(c(t))=\hat{h}_f(c) h(t)+O(1),$
where the implicit constant only depends on $f\in k(t)(z)$ and $c\in k(t)$.
    \item \label{thm3part1}
Let $\{t_n\}_{n\in\N}\subset\overline{k}$ be a non-repeating sequence of points with $\hat{h}_{f_{t_n}}(c(t_n))\to 0$.
The sequence of probability measures, weighted equally on the points in the $\Gal(\overline{k}/k)-$orbit of $t_n$, converges weakly to $\mu_{f,c,v}$ on $\PP^{1, an}_{v}$ as $n\to\infty$ for each $v\in\cM_k$. 
  \end{enumerate}
\end{theorem}
To describe two problems our Theorem \ref{vchunlikely} sheds light on, we begin with a question motivated by 
a theorem of Tate \cite{Tate} dating back to 1983.
He considered an elliptic surface $E\to X$ and a section $P: X\to E$, both defined over a number field $k$, and proved that the map $t\to \hat{h}_{E_t}(P_t)$, associating to each $t\in X(\overline{k})$ the N\'eron-Tate height of $P_t$ in the corresponding fiber $E_t$ is actually a height function on the curve $X$ corresponding to a divisor of degree equaling the geometric canonical height $\hat{h}_E(P)$ of $P\in E(k(X))$. 
More precisely, Tate showed that there exists a divisor $D=D(E, P)\in \mathrm{Pic}(X)\otimes \Q$ of degree $\hat{h}_E(P)$ such that $\hhat_{E_t}(P_t)=h_D(t)+O(1)$, as $t\in X(\overline{k})$ varies. 
Silverman \cite{Silverman:VCHI,Silverman:VCHII,Silverman:VCHIII} strengthened Tate's result to show that not only is the error term $O(1)$, but the difference behaves quite regularly. It is an open question whether the analog of Tate's result holds in the setting of arithmetic dynamics where one cannot exploit the group structure of an elliptic curve, which played a crucial role in both Tate's and Silverman's proofs. To this end Call and Silverman \cite[Theorem 4.1]{Call:Silverman} have shown that $\hat{h}_{f_t}(c(t))=\hat{h}_f(c)h(t)+o(h(t))$ as $h(t)\to\infty$. It is thus natural to ask the following question.
\begin{question}\label{vch}
\emph{Let  $K=k(X)$ be the function field of a smooth projective curve $X$ defined over a number field $k$ and let $f\in K(z)$ be a rational map of degree $d\ge 2$ and $c\in K$. Is there a divisor $D=D(f, c)\in \mathrm{Pic}(X)\otimes\R$ of degree $\hat{h}_f(c)$ such that as $t\in X(\overline{k})$ varies we have $\hhat_{f_t}(c(t))=h_D(t)+O(1)$?}
\end{question}
Ingram \cite{I13} gave an affirmative answer to Question \ref{vch} when $f$ is a polynomial in $z$. Extending Ingram's result beyond polynomial families is hard.  
To this end Ghioca, Hsia and Tucker \cite[Theorem 5.4]{GHT:preprint} proved that  the answer to Question \ref{vch} is still positive for rational maps $f(z)\in K(z)$ and starting points $c\in K$ if certain technical conditions are satisfied. 
The only other examples where the answer to Question \ref{vch} is known to be  affirmative are Latt\`es maps coming from elliptic surfaces, and when $X=\P^1$, the family $f(z)=\frac{z^d+t}{z}$ for $d\ge 2$ studied in \cite{G-M}. 
We point out here that if $X=\P^1$, then an affirmative answer to Question \ref{vch} is equivalent to the variation of heights $\hhat_{f_t}(c(t))=\hat{h}_f(c)h(t)+O(1)$, as given in part \ref{thm3part2} of Theorem \ref{vchunlikely}. 

A common feature in all the examples where the answer to Question \ref{vch} is known to be positive, is that the corresponding dynamical pairs $(f,c)$ are adelic. Hence, only finitely many local heights contribute to the difference $\hhat_{f_t}(c(t))-h_D(t)$. 
However, as seen in our Theorem \ref{genericnoadelic}, most dynamical pairs over $\P^1$ fail to be adelic. 
In part \ref{thm3part2} of Theorem \ref{vchunlikely}, we see that in the case $X=\P^1$ answering Question \ref{vch} hinges upon understanding whether all dynamical pairs are quasi-adelic, something that we believe to be true. 
Essentially, our summability condition in the definition of a quasi-adelic measure is manufactured so that the answer to Question \ref{vch} is positive.

Ghioca, Hsia and Tucker \cite{Ghioca:Hsia:Tucker, GHT:preprint} pointed out that Question \ref{vch} is related with a question in arithmetic dynamics, motivated by conjectures in arithmetic geometry in the theme of unlikely intersections; see \cite{Zannier:book} for a beautiful overview of these problems. 
Their research was motivated by a question posed by Zannier who proposed a dynamical analog of his theorem with Masser \cite{Masser:Zannier, Masser:Zannier:2, Masser:Zannier:nonsimple} in the setting of elliptic surfaces. 
The first groundbreaking result in this direction was obtained by  Baker and DeMarco \cite{BD:preperiodic} who showed that there are infinitely many $t\in\C$ such that two starting points $a, b\in \C$ are both preperiodic under the action of $f_t(z)=z^d+t$ if and only if $a^d=b^d$; in particular, $f_t(a)=f_t(b)$ and $a$ is a preperiodic point of $f_t$ if and only if $b$ is.
Their theorem was generalized in \cite{DeMarco:Mavraki, DWY:Lattes, GH:Drinfeld, Ghioca:Hsia:Tucker, GHT:preprint} to allow for more dynamical pairs and for replacing preperiodic points with points of small canonical height in the spirit of Zhang's dynamical Bogomolov conjecture. 
In this direction, the following more general conjecture was formulated; see \cite[Conjecture 1.10]{BD:polyPCF},  \cite[Conjecture 6.1]{DeMarco: heights},  \cite[Question 1.3]{Ghioca:Hsia:Tucker} and  \cite[Conjecture 2.3]{GHT:preprint}.
\begin{conjecture}[Baker-DeMarco, Ghioca-Hsia-Tucker]\label{iffconj}
\emph{Let  $K=k(X)$ be the function field of a smooth projective curve $X$ defined over a number field $k$, and consider $f\in K(z)$ and $c_1,c_2\in K$. Assume that there are infinitely many $t_n\in X(\overline{k})$ such that $\hat{h}_{f_{t_n}}(c_1(t_n))+\hat{h}_{f_{t_n}}(c_2(t_n))\to 0$ as $n\to\infty$. Then $c_1, c_2$ are \emph{coincident} along $X$, i.e., there exists $i\in\{1,2\}$ and a Zariski open subset $Y\subset X$ such that
\begin{align*}
\{t\in Y(\overline{k}):c_1(t) \text{ \bf{and} } c_2(t)\text{ is preperiodic for }f_t\}
=\{t\in Y(\overline{k}):c_i(t) \text{ is preperiodic for }f_t\}.
\end{align*}
}
\end{conjecture}
Recall that preperiodic points are points with canonical height equal to zero.  
Hence, a special case of Conjecture~\ref{iffconj} asserts that if there are infinitely many parameters $t\in X(\overline{k})$ such that both $c_1(t)$ and $c_2(t)$ are preperiodic under iteration by $f_t$, then either one of the $c_i$ is identically preperiodic for $f$, or for each $t\in Y(\overline{k})$ we have $c_1(t)$ is preperiodic under iteration by $f_t$ if and only if $c_2(t)$ is preperiodic under iteration by $f_t$. The common key ingredient in proving special cases of Conjecture \ref{iffconj}  was the adelic equidistribution theorem in \cite{BRbook, CLoir, FRL:equidistribution, Thuillier:these, Yuan}.  
However, as seen by our Theorem \ref{genericnoadelic}, most dynamical pairs $(f,c)$ on $\P^1$ are not adelic. 
In spite of this, part \ref{thm3part1} of our Theorem \ref{vchunlikely}, which is an application of Theorem \ref{arithmetic equidistribution}, implies that Conjecture \ref{iffconj} holds when $X=\P^1$ assuming that the dynamical pairs $(f,c_i)$ are quasi-adelic. 
Hence our Theorem \ref{arithmetic equidistribution} is important for applications.

Other applications of adelic equidistribution theorems include the study of the distribution of postcritically finite maps in the moduli space of rational functions, see \cite{ BD:polyPCF, DWY:Per1, Favre:Gauthier, Favre:Gauthier2, AO1, AO2, AO3}. Part \ref{thm3part1} of our Theorem \ref{vchunlikely} has various implications in this setting as well. 

\medskip

\noindent {\bf Outline of the article.} In Section \S\ref{definitions}, we introduce the notion of a quasi-adelic measure on $\PP^1$ and study some properties of this measure. In Section \S\ref{proofs of main}, we prove our main equidistribution Theorem \ref{arithmetic equidistribution} and establish an important finiteness property for the height associated with a quasi-adelic measure (see Proposition \ref{finiteness}). In Section \S\ref{a pair is rarely adelic} we prove Theorems \ref{cornoadelic} and \ref{genericnoadelic}. In Section \S\ref{quasi-adelic family}, we prove Theorem \ref{mainthm3}, thus giving  examples of dynamical pairs $(f, c)$ that are quasi-adelic, but fail to be adelic. Finally, in Section \S\ref{Svchunlikely} we prove Theorem \ref{vchunlikely}.

\medskip

\noindent {\bf Acknowledgements.}
We are indebted to Dragos Ghioca for his support and for various comments and suggestions on previous versions of this article that greatly improved its exposition.
We are also grateful to Matthew Baker, Laura DeMarco, Patrick Ingram, Xiaoguang Wang and Khoa Nguyen for many helpful comments and suggestions. We also thank the Fields Institute for its hospitality. Many of the results of this article were finalized there.  
\section{Quasi-adelic measure and some properties}\label{definitions}
In this section, we introduce the notion of quasi-adelic measure and quasi-adelic set. Further we define canonical heights associated with these measures and prove some of their properties.
 \subsection{Preliminaries and basic notation}
A {\em product formula field} is a field $k$ together with a set $\cM_k$ consisting of pairwise inequivalent non-trivial absolute values, and a unique positive integer $N_v$ associated to each element of $\cM_k$ such that the following holds.
\begin{itemize}
\item For each $\alpha \in k^*$, we have $|\alpha|_v=1$ for all but finitely many places $v\in \cM_k$, and the {\em product formula} holds
\begin{equation} \label{product formula}
\prod_{v \in \cM_k} |\alpha|_v^{N_v} \ = \ 1 \ .
\end{equation}
\end{itemize}
In what follows, we often refer to the elements of $\cM_k$ as places of $k$. 
Important examples of product formula fields include number fields and function fields of smooth projective curves. Let $\kbar$ and $\ksep$ be the algebraic and respectively separable closure of $k$. If the characteristic of $k$ is zero, then $\kbar=\ksep$. 
For each $v \in \cM_k$, let $k_v$ be the completion of $k$ with respect to $|\cdot|_v$, $\kvbar$ be an algebraic closure of $k_v$ and $\CC_v$ denote the completion of $\kvbar$. 
 We also let $\PP^{1, an}_{v}$ be the Berkovich projective line over $\CC_v$.    This is a canonically defined path-connected compact Hausdorff space containing $\PP^1(\CC_v)$ as a dense subspace. For each $v\in\cM_k$, we fix an embedding of $\kbar$ into $\C_v$. 
 We remark here that if $v$ is archimedean, we have $\CC_v \simeq \CC$ and $\PP^{1, an}_{v} \simeq \PP^1(\CC)$.

For each $v \in \cM_k$ there is a distribution-valued Laplacian operator $\Delta$ on $\PP^{1, an}_{v}$. For its definition and some examples we refer the reader to  \cite[Chapter 5]{BRbook}. 
An  important example is the Laplacian of $\log^{+}|z|_v:= \max\{\log  |z|_v, 0\}$. Note that the function
$\log^+|z|_v$, which is originally defined on $\PP^1(\CC_v)\setminus\{\infty\}$, extends naturally to a continuous real valued function defined on $\PP^{1, an}_{v} \backslash \{ \infty \}$. The Laplacian of its extension, also denoted by $\log^{+}|z|_v$, is
\begin{equation}\label{lambdav potential}
\Delta \log^+|z|_v = \delta_{\infty} - \lambda_v,
\end{equation}
where $\lambda_v$ is the uniform probability measure on the complex unit circle $\{ |z|_v = 1 \}$ when $v$ is archimedean and a point mass at the Gauss point of $\PP^{1, an}_{v}$ when $v$ is non-archimedean.

A probability measure $\mu_v$ on $\PP^{1, an}_{v}$ is said to have {\em continuous potentials} if $\mu_v - \lambda_v = \Delta g$
for some continuous function $g : \PP^{1, an}_{v} \to \RR$.  We call the funtion $g$ a potential of $\mu_v$ and note that any two potentials of $\mu_v$ differ by a constant.  If $\mu_v$ has continuous potentials,  then there is a \textbf{unique} function $g^o_{\mu_v}: \PP^{1, an}_{v} \to \RR\cup\{+\infty\}$ such that the function
    $$g(z)=\log^+|z|_v-g^o_{\mu_v}(z),$$
is a continuous potential for $\mu_v$ \textbf{and} the following normalization condition holds. If $G_{\mu_v}(x,y)$ is the function uniquely determined by $g^o_{\mu_v}$ as 
\begin{equation}
G_{\mu_v}(x,y): =\begin{cases}
g^o_{\mu_v}(x/y)+\log|y|_v, \textup{ for $x,y\in \C_v$ and $y\neq 0
$}\\
\log|x|_v-g(\infty),\textup{ for $x\in \C_v$, $x\neq 0$ and $y=0$}\\
-\infty , \textup{ for $x=y=0$},
\end{cases}
\end{equation}
where $g(\infty)=\displaystyle\lim_{z\to \infty}(\log^+|z|_v-g^o_{\mu_v}(z))$, 
then the set
\begin{equation}\label{Kv set}
M_{\mu_v}:=\{(x,y)\in \C_{v}^2 ~:~ G_{\mu_v}(x,y)\leq 0\},
\end{equation}
has capacity $\capacity(M_{\mu_v})=1$; see \cite[\S 3.3]{Baker:Rumely:equidistribution} and \cite{D:bifurcations} for the definition of homogeneous capacity. 
We call the continuous function $G_{\mu_v}$ on $\C_v^2\backslash \{(0, 0)\}$ the {\em normalized homogeneous potential} of $\mu_v$. Any other homogeneous potential of $\mu_v$ differs from $G_{\mu_v}$ by a constant. 
For example, from (\ref{lambdav potential}) we get that the normalized homogeneous potential of $\lambda_v$ is
  $$G_{\lambda_v}(x,y)=\log \|(x,y)\|_v\textup{ and } M_{\lambda_v}=\bar{D}^2(0,1)\subset \C_v^2,$$
where $\|\cdot\|_v$ is the maximum norm defined as $\|(x,y)\|_v:=\max\{|x|_v, |y|_v\}.$
An important property of a homogeneous potential $G_{\mu_v}$ is that it {\em scales logarithmically}
   $$G_{\mu_v}(\alpha x,\alpha y)=G_{\mu_v}(x, y)+\log|\alpha|_v.$$

Finally, for each probability measure $\mu_v$ on $\P^{1, an}_{v}$ with continuous potentials we have a unique {\em normalized Arakelov-Green function} $g_{\mu_v} : \PP^{1, an}_{v} \times \PP^{1, an}_{v} \to \RR \cup \{ +\infty \}$. This is characterized by the differential equation $\Delta_x g_{\mu_v}(x,y) = \delta_y - \mu_v$ and the normalization
\begin{equation} \label{normalization}
	\iint g_{\mu_v}(x,y) d\mu_v(x) d\mu_v(y) = 0.
\end{equation}
For points $(x,y)\in \C_v^2$, the normalized Arakelov-Green function $g_{\mu_v}$ is given by
\begin{equation}  \label{explicit g}
g_{\mu_v}(x,y) = -\log|\tilde{x} \wedge \tilde{y}|_v + G_{\mu_v}(\tilde{x}) + G_{\mu_v}(\tilde{y}),
\end{equation}
where $\tilde{x}=(x_1,x_2)$ and $\tilde{y}=(y_1,y_2)$ are lifts of $x$ and $y$ respectively and  $|\tilde{x}\wedge \tilde{y}|_v:=|x_1y_2-y_1x_2|_v$. Since $\PP^1(\C_v)$ is dense in $\P^{1, an}_{v}$, by continuity, we see that the extension of the function $g_{\mu_v}$ on $\PP^{1, an}_{v} \times \PP^{1, an}_{v}$ is uniquely determined by $G_{\mu_v}$.  
The fact that the capacity of the set $M_{\mu_v}$ defined in (\ref{Kv set}) is equal to $1$ guarantees that the function $g_{\mu_v}$ defined in (\ref{explicit g}) satisfies the integral formula (\ref{normalization}); see \cite{Baker:Rumely:equidistribution,BRbook} for more details.
\subsection{Quasi-adelic measure and canonical height function}\label{Quasi-adelic measure and canonical height function}
Let $\mu_v$ be a probability measure on $\P^{1, an}_{v}$ with continuous potentials. We define the {\em outer radius} and  {\em inner radius} for $\mu_v$ as
\begin{align*}   
r_{out}(\mu_v)&:=\inf \{r>0~:~ M_{\mu_v} \subset \bar{D}(0, r)\times \bar{D}(0,r)\}\\
 r_{in}(\mu_v)&:=\sup \{r>0~:~ \bar{D}(0, r)\times \bar{D}(0,r) \subset M_{\mu_v}\}.
\end{align*}
A {\em quasi-adelic measure} on $\PP^1$ with respect to a product formula field $k$ is a collection ${\mathbb \mu} = \{ \mu_v \}_{v \in \cM_k}$ of probability measures
on $\PP^{1, an}_{v}$, one for each $v \in \cM_k$, such that
\begin{itemize}
\item  $\mu_v$ has continuous potentials for each $v \in \cM_k$, and
\item $\displaystyle\prod_{v\in \cM_k} r_{in}(\mu_v)^{N_v}>0$ and  $\displaystyle\prod_{v\in \cM_k} r_{out}(\mu_v)^{N_v} <\infty$.
\end{itemize}
\begin{remark} Since  $\capacity(M_{\mu_v})=1$ and $\capacity(\bar{D}^2(0,r))=r^2$, the radii satisfy $0< r_{in}(\mu_v)\leq 1\leq  r_{out}(\mu_v)$. The measure ${\mathbb \mu} = \{ \mu_v \}_{v \in M_k}$ is {\em adelic} if we replace the  second condition by $\mu_v=\lambda_v$ or equivalently $r_{in}(\mu_v)= r_{out}(\mu_v)=1$ for all but finitely many $v\in \cM_k$; see \cite{BRbook, FRL:equidistribution}.  In other words, adelic measures satisfy $M_{\mu_v}=\bar{D}^2(0,1)$ for all but finitely many places $v\in\cM_k$.
\end{remark}
If $\rho,\rho'$ are probability measures on $\PP^{1, an}_{v}$, we define the $\mu_v$-energy of $\rho$ and $\rho'$ as 
$$( \rho, \rho' )_{\mu_v} := \frac{1}{2} \iint_{\PP^{1, an}_{v} \times \PP^{1, an}_{v} \backslash {\rm Diag}} g_{\mu_v}(x,y) d\rho(x) d\rho'(y).$$
 The {\em $\mu_v$-energy of $\rho$} is defined as $I_{\mu_v}(\rho) := (\rho,\rho)_{\mu_v}$.

We can now define the height associated with a quasi-adelic measure.
Let $S\subset \PP^1(\ksep)$ be a finite Gal$(\ksep/k)$-invariant set with cardinality $|S|>1$. Let $\tilde{S}$  be a Gal$(\ksep/k)$-invariant set consisting of lifts $\tilde{x}$ of elements $x\in S$; in particular $|\tilde{S}|=|S|$. 
We denote by $[S]_v$ the discrete probability measure on $\PP^{1, an}_{v}$ supported
equally on all elements of $S$. The {\em canonical height} of $S$ associated to a quasi-adelic measure ${\mathbb \mu} = \{ \mu_v \}_{v \in M_k}$ is a number given by
\begin{align}\label{height definition}
\begin{split}
\hhat_{\mu}(S) :&= \frac{|S|}{|S|-1}\sum_{v\in \cM_k} N_v \cdot ([S]_v, [S]_v)_{\mu_v}  \\
	&=  \frac{|S|}{|S|-1}  \sum_{v\in \cM_k} \frac{N_v }{2|S|^2}\sum_{x,y\in S, x\not=y}  g_{\mu_v}(x,y) \\	
	&=  \sum_{ x,y\in S, x\not=y } ~  \sum_{v\in \cM_k}  \frac{N_v \cdot ( -\log|\tilde{x} \wedge \tilde{y}|_v + G_{\mu_v}(\tilde{x}) + G_{\mu_v}(\tilde{y}) )}{2|S|(|S|-1)}\\
	&=  \frac{1}{2|S|(|S|-1)} \sum_{\tilde{x},\tilde{y}\in \tilde{S}, x\not=y}  \; \sum_{v\in \cM_k} N_v \cdot (G_{\mu_v}(\tilde{x}) + G_{\mu_v}(\tilde{y}) )  \textup{, by (\ref{product formula})}\\
     &=  \frac{1}{|S|}  \cdot \sum_{\tilde{x}\in \tilde{S}}\sum_{v\in \cM_k} N_v\cdot  G_{\mu_v}(\tilde{x}),   
\end{split}
\end{align}
Here the constants $N_v$ are the same as those appearing in the product formula.
Therefore, we have
\begin{align}\label{heightq}
\hhat_{\mu}(S)= \frac{1}{|S|}  \cdot \sum_{\tilde{x}\in \tilde{S}}\sum_{v\in \cM_k} N_v \cdot G_{\mu_v}(\tilde{x}).
\end{align}
If $x\in \ksep$ we may take $S=\Gal(\ksep/ k)\cdot x$ and use \eqref{heightq} to define the canonical height of $x$ as
   $$\hhat_{\mu}(x) :=\hhat_{\mu}(S).$$
\begin{remark} The canonical height we defined is slightly different from the one appeared in \cite{BRbook, FRL:equidistribution}, but agrees with the one in \cite{DWY:Lattes, DWY:Per1}. The factor of $\frac{|S|}{|S|-1}$ is included here to allow for a better comparison of this measure-theoretic height with the Call-Silverman height; see Proposition \ref{samesmallpoints}. Equation \eqref{heightq} allows us to extend the definition of our height to the case $|S|=1$.
\end{remark}
\subsection{ Quasi-adelic set} In this section we introduce the notion of a quasi-adelic set. This has a geometric interpretation; hence in many applications it is easier to manipulate than a quasi-adelic measure.
 Analogous to the notion of the {\em homogeneous filled Julia set} in \cite[\S 3.2]{Baker:Rumely:equidistribution}, we define a {\em homogeneous set with continuous potential} as
   $$M_v:=\{(x,y)\in \C_v^2~ : ~ G_{M_v}(x,y)\leq 0\},$$
where $G_{M_v}$ is a continuous homogeneous potential for a  probability measure  on $\P^{1, an}_{v}$. We denote this measure, corresponding to $M_v$, by $\mu_{M_v}$. 
When $v$ is an archimedean place, $M_v$ having continuous potential is equivalent with saying that $M_v\subset \C_v^2\simeq \C^2$ is a compact, circled and pseudoconvex set, or that $G_{M_v}$ is a continuous and plurisubharmonic function satisfying
\begin{enumerate}
\item	 $G_{M_v}(\alpha z) =G_{M_v}(z) + \log|\alpha|_v$ for all $\alpha \in \C_v$, and
\item $G_{M_v}(z) = \log\|z\|_v + O(1)$;
\end{enumerate}
see \cite{D:bifurcations}.
We point out here that there are many homogeneous sets with continuous potential.
 If $F_n: \C_v^2\to \C_v^2$ is a sequence of homogeneous polynomials with  $\deg (F_n)\geq 1$  such that the sequence of functions $\{\frac{\log \| F_n\|_v}{\deg (F_n)}\}_{n\geq 1}$ converges uniformly to $G_v$ on $\C_v^2\backslash \{(0, 0)\}$, then $G_v$ is a homogeneous potential for some probability measure with continuous potentials on $\P^{1, an}_{v}$; see \cite[\S 3]{Baker:Rumely:equidistribution}.
 Hence $M_v=\{(x,y) \in \C_v^2: G_v(x,y)\leq 0\}$ is a homogeneous set with continuous potential. As seen in \cite[$\S$2]{DWY:Lattes}, its capacity can be computed by the following limit
   $$\capacity(M_v)=\lim_{n\to \infty} |\Res(F_n)|_v^{-\frac{1}{\deg(F_n)^2}}.$$
Analogous to the definition of the radii of $\mu_v$, we define the {\em outer} and {\em inner} radii of $M_v$ as
   $$r_{out}(M_v):=\inf \{r>0~:~ M_{v} \subset \bar{D}(0, r)\times \bar{D}(0,r)\}$$
    $$r_{in}(M_v):=\sup \{r>0~:~ \bar{D}(0, r)\times \bar{D}(0,r) \subset M_{v} \}.$$
A product $\displaystyle\prod_{v\in \cM_k}r_v^{N_v}$, with $r_v> 0$ for each $v\in \cM_k$,  \emph{converges strongly} if
   $$\sum_{v\in \cM_k}N_v\cdot |\log r_v|<\infty.$$
We define a {\em quasi-adelic set}  (with respect to a product formula field $k$) to be a collection  ${\mathbb M} = \{ M_v \}_{v \in \cM_k}$ of sets such that the following hold.
\begin{itemize}
\item For each $v\in\cM_k$ the set $M_v$ is a homogeneous set with continuous potential.
\item The products $\displaystyle\prod_{v\in \cM_k} r_{out}(M_v)^{N_v}$ and $\displaystyle\prod_{v\in \cM_k} r_{in}(M_v)^{N_v}$ converge strongly.
\end{itemize}
Note that there is a unique probability measure $\mu_{M_v}$ with continuous potential associated to a homogeneous set $M_v$ with continuous potential. Hence a quasi-adelic set ${\mathbb M} = \{ M_v \}_{v \in M_k}$ gives a measure ${\mathbb \mu_{\mathbb M}} = \{ \mu_{M_v} \}_{v \in M_k}$ on $\PP^1$. In the next theorem we will see that this measure is also quasi-adelic.
\begin{theorem}\label{quasi-adelic set}
Let $k$ be a product field and $\mathbb M=\{M_v\}_{v\in  \cM_k}$ be a collection of homogeneous sets with continuous potential. Then we have:
\begin{itemize}
\item If the set $\mathbb M=\{M_v\}_{v\in  \cM_k}$ is quasi-adelic, then the corresponding measure ${\mathbb \mu_{\mathbb M}} = \{ \mu_{M_v} \}_{v \in M_k}$ is quasi-adelic.
\item Suppose that for each $v\in \cM_k$, there are positive constants $r'_v, r_v$ such that $\bar{D}^2(0,r'_v)\subset M_v\subset \bar{D}^2(0,r_v)$
and the products $\prod_{v\in \cM_k}{r'_v}^{N_v}$, $\prod_{v\in \cM_k}{r_v}^{N_v}$ converge strongly. Then the set $\mathbb M=\{M_v\}_{v\in  \cM_k}$ is quasi-adelic. Moreover, the product $\prod_{v\in \cM_k}\capacity(M_v)^{N_v}$ converges strongly and for any $\Gal$($\ksep/k$)-invariant $S\subset \P^1(\ksep)$ we have
 \begin{equation}\label{height cap}
 \hhat_{\mu_{\mathbb M}}(S)= \frac{1}{|S|}  \cdot \sum_{\tilde{x}\in\tilde{ S}}\sum_{v\in \cM_k} N_v \cdot G_{M_v}(\tilde{x})+\frac{1}{2}\log \prod_{v\in \cM_k}\capacity(M_v)^{N_v}.
 \end{equation}
\end{itemize}
\end{theorem}

\subsection{Some properties}
A number field $k$ is naturally equipped with a set of inequivalent absolute values $\cM_k$ and positive integers $\{N_v\}_{v\in \cM_k}$ making it a product formula field.  For any $x\in k$, the {\em logarithmic Weil height} of $x$ is defined as follows
\begin{align}\label{logweilheight} 
  h(x):=\frac{1}{[k:\Q]}\sum_{v\in \cM_k} N_v\cdot \log^+|x|_v.
\end{align}
It is well defined and does not depend on the embedding of a number field $k\hookrightarrow\Qbar$. An important property for an adelic measure is that its canonical height differs from a multiple of the Weil height by a bounded amount. The following proposition establishes a similar result for quasi-adelic measures.

\begin{prop} \label{equivalent heights} Let $k$ be a number field. Suppose  $\mu=\{\mu_{v}\}_{v\in \cM_k}$ is a quasi-adelic measure. Then the canonical height $\hhat_{\mu}$ is bounded by the logarithmic Weil height $h$ on $\PP^1(\kbar)$ as
     $$\log \prod_{v\in \cM_k} r_{in}(\mu_v)^{N_v}\leq [k:\Q]h(x)-\hhat_\mathbb{\mu}(x)\leq \log \prod_{v\in \cM_k} r_{out}(\mu_v)^{N_v},$$
for all $x\in \kbar$.
\end{prop}
\begin{proof}  Let $x\in\kbar$ and write $S=\Gal(\kbar/ k)\cdot x$ to for its Galois orbit and  $\tilde{x}\in \kbar^2$ for a lift of $x\in \kbar$. The definition of outer and inner radii gives
   $$\log r_{in}(\mu_v) \leq\frac{1}{|S|}\sum_{y\in S} (\log\|\tilde{y}\|_v-G_{\mu_v}(\tilde{y}))\leq \log r_{out}(\mu_v).$$
Consequently,
   $$\sum_{v\in \cM_k} N_v\log r_{in}(\mu_v) \leq\sum_{v\in \cM_k} \frac{N_v}{|S|}\sum_{y\in S} (\log\|\tilde{y}\|_v-G_{\mu_v}(\tilde{y}))\leq \sum_{v\in \cM_k} N_v\log r_{out}(\mu_v).$$
Using (\ref{heightq}) this inequality can be rewritten as
$$\log \prod_{v\in \cM_k} r_{in}(\mu_v)^{N_v}\leq [k: \Q]h(x)-\hhat_\mathbb{\mu}(x)\leq \log \prod_{v\in \cM_k} r_{out}(\mu_v)^{N_v}.$$\end{proof}
Finally, we mention that like the tensor product of adelic metrized line bundles is again an adelic metrized line bundle, it is easy to see that certain linear combinations of quasi-adelic measures are quasi-adelic measures. For example, the average of two quasi-adelic measures is a quasi-adelic measure. 

\section{Equidistribution of small points}\label{proofs of main}
In this section, we  prove Theorems \ref{arithmetic equidistribution} and \ref{quasi-adelic set}. Moreover, we prove an important finiteness property for the height associated to a quasi-adelic measure. 

\subsection{Proof of Theorem \ref{arithmetic equidistribution}} By assumption,  $\{S_n\}_{n\geq 1}$ is a sequence of subsets of $\PP^1(\ksep)$ which are Gal$(\ksep/k)$-invariant and the cardinality $|S_n|$ tends to infinity. For each $v\in \cM_k$, the $\mu_v$-energy of the probability measure $[S_n]_v$ on $\P^{1, an}_{v}$ is given by
    \begin{align*}([S_n]_v, [S_n]_v)_{\mu_v}&=\frac{1}{2|S_n|   ^2}\sum_{x\neq y\in S_n} g_{\mu_v}(x,y)\\
      &=\frac{1}{2|S_n|   ^2}\sum_{x\neq y\in S_n} \left(-\log|\tilde{x} \wedge \tilde{y}|_v + G_{\mu_v}(\tilde{x}) + G_{\mu_v}(\tilde{y})\right),
\end{align*}
where $\tilde{x}$ and $\tilde{y}$ are lifts of $x$ and $y$ respectively. We begin by proving the following lemma.

\begin{lemma}\label{energy limit}
Let $\{S_n\}_{n\geq 0}$ and $\mu=\{\mu_v\}_{v\in \cM_k}$ be as in Theorem \ref{arithmetic equidistribution}. For each place $v\in \cM_k$, we have
   $$\lim_{n\to \infty} ([S_n]_v, [S_n]_v)_{\mu_v}=0.$$
\end{lemma}
\proof 
First we will show that for each $v\in \cM_k$, we have
\begin{equation}\label{inf limit}
\lim_{n\to \infty}\inf  ([S_n]_v, [S_n]_v)_{\mu_v}\geq 0.
\end{equation}
For this we follow the proof of \cite[Lemma 3.17]{Baker:Rumely:equidistribution}. From the definition of homogenous capacity we have
\begin{equation}\label{22222}\lim_{n\to \infty} \inf \inf_{\tilde x_1, \tilde x_2, \ldots, \tilde x_n \in M_{\mu_v}}\frac{1}{n(n-1)}\sum_{i\neq j}-\log|\tilde{x_i}\wedge \tilde{x}_j|_v\geq -\log \capacity(M_{\mu_v})=0.
\end{equation}
   Let $\epsilon>0$ be an arbitrary small number. Since $\{|\alpha |_v: \alpha \in \C_v\}$ is dense in $\R_{\ge 0}$, we can choose lifts of $x, y\in S_n$ , denoted by $\tilde{x}, \tilde{y} \in M_{\mu_v}$, such that
    $$-\epsilon<G_{\mu_v}(\tilde{x})\leq 0\text{ and }-\epsilon<G_{\mu_v}(\tilde{y})\leq 0.$$
Now the definition of the energy function and (\ref{22222}) yield
     \begin{align*}\lim_{n\to \infty}\inf ([S_n]_v, [S_n]_v)_{\mu_v} &=\lim_{n\to \infty}\inf\frac{1}{2|S_n|   ^2}\sum_{x\neq y\in S_n} \left(-\log|\tilde{x} \wedge \tilde{y}|_v + G_{\mu_v}(\tilde{x}) + G_{\mu_v}(\tilde{y})\right)\\
     &\geq -\epsilon+\lim_{n\to \infty}\inf\left(-\frac{1}{2|S_n|   ^2}\sum_{x\neq y\in S_n} \log|\tilde{x} \wedge \tilde{y}|_v\right) \geq -\epsilon.
\end{align*}
Shrinking $\epsilon$ to zero, \eqref{inf limit} follows. 
It remains to prove that for each $v\in \cM_k$ we have
   $$\lim_{n\to \infty}\sup ([S_n]_v, [S_n]_v)_{\mu_v}\leq 0.$$
We assume that this inequality fails to end in a contradiction. Then there exist an $\epsilon>0$, $v_o\in \cM_k$ and a sequence of strictly increasing integers $\{n_j\}_{j\geq 1}$ such that
\begin{equation}\label{big limit} N_{v_o}\cdot ([S_{n_j}]_{v_o}, [S_{n_j}]_{v_o})_{\mu_{v_o}}>3\epsilon, \textup{ for all $j\geq 1$}.
\end{equation}
For each $v\in \cM_k$ and any $\delta>0$, we lift $x,y\in S_n$ to $\tilde{x}, \tilde{y}\in M_{\mu_v}$ such that
      $$G_{\mu_v}(\tilde{x})\geq -\delta, ~ G_{\mu_v}(\tilde{y})\geq -\delta.$$
Since we know that $M_{\mu_v}$ is bounded by the polydisc with outer radius $r_{out}(\mu_v)$ we have $\|\tilde{x}\|_v, \|\tilde{y}\|_v\leq r_{out}(\mu_v)$. Hence for a non-archimedean place $v\in \cM_k$ we have
  $$\log|\tilde{x} \wedge \tilde{y}|_v \leq \log r_{out}(\mu_v)^2,$$
  and
 \begin{align}
\begin{split}
([S_n]_v, [S_n]_v)_{\mu_v}&=\frac{1}{2|S_n|   ^2}\sum_{x\neq y\in S_n} \left(-\log|\tilde{x} \wedge \tilde{y}|_v + G_{\mu_v}(\tilde{x}) + G_{\mu_v}(\tilde{y})\right)\\
      &\geq \frac{1-|S_n|}{|S_n|}\delta- \sum_{x\neq y\in S_n}\frac{ \log r_{out}(\mu_v)}{| S_n|^2}=\frac{1-|S_n|}{|S_n|}\left(\delta+\log r_{out}(\mu_v)\right).
\end{split}
\end{align}
Shrinking $\delta$ to zero, we get that for all non-archimedean places $v\in\cM_k$:
\begin{equation}\label{1111}([S_n]_v, [S_n]_v)_{\mu_v}\geq \frac{1-|S_n|}{|S_n|}\log r_{out}(\mu_v).
\end{equation}
It is well known that for a product formula field $k$ there are only finitely many archimedean places; see \cite[Chapter 12, Theorem 3]{Ar}. Since $\mu=\{\mu_v\}_{v\in \cM_k}$ is quasi-adelic, the product $\prod_{v\in \cM_k} r_{out}(\mu_v)^{N_v}$ converges. Hence we can choose a set $\cM_k'\subset \cM_k$ such that:
\begin{itemize}
\item $\cM_k'':=\cM_k\backslash( \cM_k'\cup \{v_o\})$ has only finitely many places.
\item all places in $\cM_k'$ are non-archimedean and $v_o\not \in \cM_k'$.
\item $\displaystyle\sum_{v\in \cM_k'}N_v\cdot \log r_{out}(\mu_v)\leq \epsilon$.
\end{itemize}
Now the definition of the canonical height $\hat h_\mu$ gives
 \begin{align*} N_{v_o}([S_{n_j}]_{v_o}, [S_{n_j}]_{v_o})_{\mu_{v_o}}&=\hat h_{\mu}(S_{n_j})-\sum_{v\in \cM_k\backslash \{v_o\}} N_v \cdot([S_{n_j}]_v, [S_{n_j}]_v)_{\mu_v}\\
  &=\hat h_{\mu}(S_{n_j})-\sum_{v\in \cM_k''} N_v ([S_{n_j}]_v, [S_{n_j}]_v)_{\mu_v}-\sum_{v\in \cM_k'} N_v  ([S_{n_j}]_v, [S_{n_j}]_v)_{\mu_v}.
\end{align*}
This in turn, upon using \eqref{1111} implies
\begin{align*} N_{v_o}([S_{n_j}]_{v_o}, [S_{n_j}]_{v_o})_{\mu_{v_o}}&\leq \hat h_{\mu}(S_{n_j})-\sum_{v\in \cM_k''} N_v \cdot ([S_{n_j}]_v, [S_{n_j}]_v)_{\mu_v}+\frac{(| S_{n_j}|-1) \epsilon}{| S_{n_j}|}\\
  &\leq \hat h_{\mu}(S_{n_j})-\sum_{v\in \cM_k''} N_v \cdot ([S_{n_j}]_v, [S_{n_j}]_v)_{\mu_v}+ \epsilon.
\end{align*}
Since the set $\cM_k''$ contains only finitely many places and the height of $S_{n_j}$ tends to zero,  by (\ref{inf limit}), taking the superior limit in the above inequality yields
\begin{align}\label{delta small}
\lim_{j\to\infty} \sup N_{v_o}([S_{n_j}]_{v_o}, [S_{n_j}]_{v_o})_{\mu_{v_o}}&\leq \lim_{j\to \infty}\sup \left(\hat h_{\mu}(S_{n_j})-\sum_{v\in \cM_k''} N_v ([S_{n_j}]_v, [S_{n_j}]_v)_{\mu_v}+ \epsilon\right)\nonumber\\
 &\leq  \lim_{j\to \infty}\sup \hat h_{\mu}(S_{n_j})+\epsilon\leq \epsilon,
\end{align}
which contradicts with our assumption in (\ref{big limit}). This finishes the proof of the lemma. \qed

The set of probability measures on $\PP^{1, an}_{v}$ is compact in the weak topology. Hence, to show that $[S_n]_v$ converges weakly to $\mu_v$ as $n\to \infty$, it suffices to show that any convergent subsequence of  $[S_n]_v$ converges to $\mu_v$. Without loss of generality, we assume that $[S_n]_v$ converges to some $\nu_v$,
   $$\lim_{n\to \infty}[S_n]_v=\nu_v.$$
By Lemma \ref{energy limit}, the $\mu_v$-energy $I_{\mu_v}(\nu_v)$ of $\nu_v$ satisfies
\begin{align*} 0&=\lim_{n\to \infty}([S_{n}]_v, [S_{n}]_v)_{\mu_v}=\lim_{n\to\infty} \frac{1}{2} \iint_{\PP^{1, an}_{v} \times \PP^{1, an}_{v} \backslash {\rm Diag}} g_{\mu_v}(x,y) d[S_n]_v(x)d[S_n]_v(y)\\
    &\geq \frac{1}{2} \iint_{\PP^{1, an}_{v} \times \PP^{1, an}_{v}} g_{\mu_v}(x,y)  d\nu_v(x)d\nu_v(y) \textup{, by \cite[Lemma 3.26]{Baker:Rumely:equidistribution}}\\
    &=I_{\mu_v}(\nu_v).
\end{align*}
Since by \cite[Theorem 3.25]{Baker:Rumely:equidistribution}, $\mu_v$ is the unique probability measure on $\P^{1, an}_{v}$ minimizing the  $\mu_v$-energy function $I_{\mu_v}(\cdot)$ and $I_{\mu_v}(\mu_v)=0\geq I_{\mu_v}(\nu_v)$, we get that $I_{\mu_v}(\mu_v)=I_{\mu_v}(\nu_v)$ and $\nu_v=\mu_v$. This finishes the proof of Theorem \ref{arithmetic equidistribution}.
\qed
\subsection{Proof of Theorem \ref{quasi-adelic set}} Firstly, we show that $\mathbb M=\{M_v\}_{v\in  \cM_k}$ is quasi-adelic implies that ${\mathbb \mu} = \{ \mu_{M_v} \}_{v \in M_k}$ is quasi-adelic. Assume that $\mathbb M=\{M_v\}_{v\in  \cM_k}$ is a quasi-adelic set. For any $r>0$, let
   $$rM_v:=\{(\alpha x, \alpha y)~:~ (x,y)\in M_v, \alpha \in \C_v \textup{ with  } |\alpha|_v\leq r\}.$$
From the definition of the capacity we have $\capacity(rM_v)=r^2\capacity(M_v)$.
Since $G_{M_v}$ is a homogeneous potential for $\mu_{M_v}$, the normalized homogeneous potential $G_{\mu_{M_v}}$ is given by
\begin{equation}\label{potential relation}
G_{\mu_{M_v}}(x,y)=G_{M_v}(x,y)+\frac{1}{2}\log \capacity(M_v),
\end{equation}
and $M_{\mu_{M_v}}=\frac{1}{\sqrt{\capacity(M_v)}}M_v$. As a consequence,
\begin{equation}\label{radii capacity}r_{in}(\mu_{M_v})=\frac{r_{in}(M_v)}{\sqrt{\capacity(M_v)}}\text{ and } r_{out}(\mu_{M_v})=\frac{r_{out}(M_v)}{\sqrt{\capacity(M_v)}}.
\end{equation}
Moreover, as $\capacity(\bar{D}^2(0, r))=r^2$, we have
\begin{equation}\label{capacity bound}
 r_{in}(M_v)\leq \sqrt{\capacity(M_v)}\leq r_{out}(M_v).
\end{equation}
Then by (\ref{radii capacity}) we have
     $$\frac{r_{in}(M_v)}{r_{out}(M_v)}\leq r_{in}(\mu_{M_v})\leq 1\leq  r_{out}(\mu_{M_v})\leq \frac{r_{out}(M_v)}{r_{in}(M_v)}.$$
As $\mathbb M=\{M_v\}_{v\in  \cM_k}$ is quasi-adelic, the products of inner and outer radii converge strongly. Then the above inequalities imply that the products of the inner and outer radii of $\mu=\{\mu_{M_v}\}_{v\in \cM_k}$ converge, that is $\mu=\{\mu_{M_v}\}_{v\in \cM_k}$ is quasi-adelic.

Now assume that $\bar{D}^2(0,r'_v)\subset M_v\subset \bar{D}^2(0,r_v)$
and $\prod_{v\in \cM_k}{r'_v}^{N_v}$, $\prod_{v\in \cM_k}{r_v}^{N_v}$ converge strongly. Then the products $\prod_{v\in \cM_k} r_{out}(M_v)^{N_v}$ and $\prod_{v\in \cM_k} r_{in}(M_v)^{N_v}$ converge strongly, since
  $$r'_v\leq r_{in}(M_v)\leq r_{out}(M_v)\leq r_v.$$
Hence $\mathbb M=\{M_v\}_{v\in  \cM_k}$ is quasi-adelic. Moreover, by (\ref{capacity bound}), the product of the capacities converges strongly. Then  the last formula for the canonical height is clear from (\ref{height definition}) and (\ref{potential relation}). We finish the proof of Theorem \ref{quasi-adelic set}.\qed

\subsection{A finiteness property}The following proposition will be useful in the last section in proving the equidistribution of parameters $t$ with small height with respect to $\hhat_{f_{c(t)}}(c(t))$. 
\begin{prop}\label{finiteness}
Let $k$ be a product formula field. Suppose $\mathbb{\mu}$ is a quasi-adelic measure. Then for any $\delta>0$ there are at most finitely many $x\in \ksep$ with
   $$\hhat_{\mu}(x)<-\delta.$$
\end{prop}
\proof Assume to the contrary that there are infinitely many $x_i\in \ksep$ with $\hhat_{\mu}(x_i)<-\delta$. Let $S_n=\cup_{i=1}^n \Gal(\ksep/k)\cdot x_i$. Then $S_n$ is $\Gal(\ksep/k)$-invariant and $|S_n|\to \infty$. 
By \eqref{height definition} we get $\hhat_{\mu}(S_n)\leq -\delta$. Moreover, from (\ref{inf limit}) we see that for each $v\in \cM_k$
   $$\lim_{n\to \infty} \inf ([S_n]_v, [S_n]_v)_{\mu_v}\geq 0.$$
Replacing $\hhat_{\mu} (S_{n_j})$ in \eqref{delta small} by $\hhat_\mu(S_n)$ and letting $\epsilon$ tend to zero, we get
   $$\lim_{n\to \infty} \sup ([S_n]_v, [S_n]_v)_{\mu_v}\leq -\delta, $$
 which is a contradiction. \qed
\section{A dynamical pair $(f, c)$ on $\PP^1$ is rarely adelic}\label{a pair is rarely adelic}
In this section we aim to prove Theorems \ref{cornoadelic} and \ref{genericnoadelic}. We first introduce some notations and terminologies.
\subsection{A dynamical pair on $\P^1$}
Let $k$ be a product formula field and let $K=k(t)$. Recall that a dynamical pair on $\P^1$ is a pair $(f, c)\in K(z)\times K$ with $d=\deg_z f\geq 2$. 
We say that the pair $(f, c)$ is {\em isotrivial} if there is a family of M\"obius transformations $M_t(z)\in \overline K(z)$ such that both $M_t\circ f_t\circ M_t^{-1}(z)$ and $M_t(c(t))$ are independent of $t$. Moreover, we say that $(f, c)$ is {\em preperiodic} if the starting point $c\in K$ is preperiodic under $f\in K(z)$, that is if  there are integers $m>n\geq 0$ with $f^n(c)=f^m(c)\in K$. 

Recall further that we have a {\em canonical height} function associated to $f\in K(z)$, denoted by $\hat{h}_f:\P^1(\overline{K})\to \P^1(\overline{K})$, determined uniquely by the properties $\hat{h}_f(f(c))=d\cdot \hat{h}_f(c)$ and $\hat{h}_f(c)=h(c)+O(1)$. Alternatively, for $c\in K$ we can compute the canonical height as:
\begin{align*}
\hhat_f(c):=\lim_{n\to \infty}\frac{\deg_t f^n(c)}{d^n}.
\end{align*}
We note that $\hhat_f(c)\geq 0$ and equality holds if and only if $(f, c)$ is either isotrivial or preperiodic; see \cite{Baker, DeMarco: heights}.
\subsection{Homogenization}
Let $(f, c)\in K(z)\times K$ be a dynamical pair with degree $d=\deg_z f\geq 2$. 
In what follows we choose lifts of $f$ and $c$ on $k^2$ with homogenous parameters $(t_1,t_2)$ as follows. For a lift of $f$ we write
\begin{align*}
F_{t_1,t_2}(z,w)=(P_{t_1,t_2}(z,w),Q_{t_1,t_2}(z,w)),
\end{align*}
where $P_{t_1,t_2}, Q_{t_1,t_2}$ are homogeneous polynomials in $(z,w)$ of degree $d$, with coefficients homogeneous polynomials in $(t_1,t_2)$ of the same degree that are relatively prime. 
For a lift of $c$ we let
$$C(t_1,t_2)=(A(t_1,t_2),B(t_1,t_2)),$$
where $A$ and $B$ are homogeneous polynomials in $k[t_1,t_2]$ of the same degree and have no common linear factor in $\overline{k}[t_1,t_2]$. Moreover, we write the $n$-th iterate of $C(t_1,t_2)$ under $F_{t_1, t_2}$ as
   $$F^n_{t_1,t_2}(C(t_1,t_2))=:(A_{C,n}(t_1,t_2),B_{C,n}(t_1,t_2)).$$
    In other words, we have
\begin{align*}
A_{C,0}(t_1,t_2)=A(t_1,t_2)\text{ and }B_{C,0}(t_1,t_2)=B(t_1,t_2),
\end{align*}
while for all $n\ge 0$ we have
\begin{align}\label{recursionnot}
\begin{split}
A_{C,n+1}(t_1,t_2)&=P_{t_1,t_2}(A_{C,n}(t_1,t_2),B_{C,n}(t_1,t_2))\\
B_{C,n+1}(t_1,t_2)&=Q_{t_1,t_2}(A_{C,n}(t_1,t_2),B_{C,n}(t_1,t_2)).
\end{split}
\end{align}
Note that $f_t^n(c(t))=\frac{A_{C,n}(t,1)}{B_{C,n}(t,1)}$.

In what follows we identify $\P^1(\overline{k})$ with $\A^1(\overline{k})\cup\{\infty\}$. We say that $f\in K(z)$ \emph{degenerates} at $t\in\A^1(\overline{k})\cup\{\infty\}$ if $\deg_z(f_t)<d$ and write \emph{$\Sing(f)$} for the set of degenerating parameters. More specifically, we denote the resultant of the homogeneous polynomials $P_{t_1, t_2}, Q_{t_1, t_2}$ in $(z,w)$ as
$$\Res(F_{t_1,t_2}):=\Res_{(z,w)}(P_{t_1,t_2},Q_{t_1,t_2})\in k[t_1,t_2]\setminus\{0\},$$
and note that $\Res(F_{t_1,t_2})\in k[t_1,t_2]$ is a homogeneous polynomial and that $f$ degenerates exactly at $t=[t_1: t_2]\in \PP^1(\bar k)$ with $\Res(F_{t_1, t_2})=0$. Thus $\Sing(f)\subset \PP^1(\kbar)$ is given by
\begin{align*}
\Sing(f)=\{\alpha=t_1/t_2\in \P^1(\overline{k})~:~\Res(F_{t_1, t_2})=0\}.
\end{align*}
We also work with a lift of $f^n(c)$ defined by coprime homogeneous polynomials.
To write the greatest common divisor of $F_{t_1,t_2}^n(C)$, for each $\alpha\in \PP^1(\kbar)$ we let 
 \[ u_{\alpha}(t_1,t_2) = \left\{
\begin{array}{ll}
t_1-\alpha t_2&  \mbox{ if $\alpha\in\A^1(\kbar)$}\\
t_2 & \mbox{ if $\alpha=\infty$.}\\
\end{array} \right. \]
Moreover we let $m_{C,n}(\alpha)$ be the maximal integer $m\in\N$ with $u_{\alpha}^m|A_{C,n}$ and $u_{\alpha}^m|B_{C,n}$. Then 
$${\rm gcd} (F^n_{t_1, t_2}(C))=\prod_{\substack{\alpha\in\Sing(f)}}u_{\alpha}(t_1,t_2)^{m_{C,n}(\alpha)}.$$
We point out here that for each $\alpha\in \Sing(f)$, the sequence  $\{m_{C,n}(\alpha)/d^n\}_{n\in\N}$ converges as $n\to\infty$.
We associate with lifts $F_{t_1,t_2}$ and $C$ of $f$ and $c$ respectively, a lift of $f^n(c)$ given by coprime homogeneous polynomials in the variables $(t_1,t_2)$, defined by
$$F_{C,n}(t_1,t_2):=F^n_{t_1,t_2}(C)/{\rm gcd} (F^n_{t_1, t_2}(C)).$$
Note that $\deg F_{C,n}(t_1, t_2)=\deg_t f_t^n(c(t))$. 

Next we introduce measures associated with each dynamical pair.

\subsection{Bifurcation measure} \label{bifurcation measure}
Let $v\in\cM_k$ be an archimedean place. Then $\|\cdot\|_v$ is a Euclidean norm, $\C_v\simeq\C$ and $\P^{1, an}_{v}\simeq \PP^1(\C)$. 
Suppose that $f\in K(z)$ does not degenerate at $t_0\in \PP^1(\C)$. The dynamical pair $(f, c)$ is {\em stable} at $t_0$ when the sequence of holomorphic maps $\{t\mapsto f_t^n(c(t))\}$
forms a normal family in a neighborhood of $t_0$. The failure of normality determines a positive measure on the parameter space, the {\em bifurcation measure}. 
Let $F$ and $C$ be lifts of $f$ and $c$ respectively. The bifurcation measure on $\P^1(\C)\backslash \Sing(f)$, denoted by $\mu_c$, is constructed as
$$\mu_c:=dd^c\left(\lim_{n\to \infty}\frac{1}{d^n}\log\|F_t^n(C(t))\|_v\right),$$
and is independent of our choice of lifts.
Its support, $\supp(\mu_c)$, is exactly the set of parameters $t$ at which $(f, c)$ is unstable. Bifurcation is important in dynamics. The family $f_t$ is {\em stable} at $t_0$ if the Julia set moves holomorphically for a small perturbation of $t$ at $t_0$, or equivalently if the dynamical pair $(f, c)$ (upon passing to a finite branched cover of $\P^1$) is stable at $t_0$ for each critical point $c$; see \cite{Mc2, MSS83}.  We refer the reader to \cite{D:current, D:bifurcations, DF08} for more details. 
\subsection{Measure associated to a dynamical pair}\label{dynamicalmeasure}
Let $k$ be a product formula field and $K=k(t)$ as before. For each $v\in\cM_k$, we let
    $$G_{F,C,v}(t_1,t_2):=\lim_{n\to \infty} \frac{\log\|F_{C,n}(t_1, t_2)\|_v}{\deg F_{C,n}}.$$
This sequence converges locally uniformly on $\C_v^{2}\setminus\{(t_1, t_2): t_1/t_2\in \Sing(f) \textup{ or } t_1{=}t_2{=}0\}$; see \cite{Branner:Hubbard:1}. For each $v\in \cM_k$ and if $\deg F_{C, n}$ is not zero, there is a probability measure $\mu_{n, v}$ on $\P^{1, an}_{v}$ associated to $F_{C, n}$, which is independent of the choice of lifts for $f$ and $c$.
\begin{proposition}\label{measure for pair} Let $(f,c)\in K(z)\times K$ be a non-isotrivial and non-preperiodic dynamical pair. For each $v\in \cM_k$, the sequence of measures $\mu_{n,v}$ on $\P^{1, an}_{v}$ converges weakly to a unique probability measure $\mu_{f, c, v}$ as $n\to \infty$. Moreover, $\mu_{f, c, v}$ has continuous potentials if and only if $G_{F,C,v}$ extends continuously on $\C_v^2\backslash\{(0, 0)\}$. 
\end{proposition}
\proof First we assume that $v\in \cM_k$ is an archimedean place, so that $\C_v\simeq\C$ and  $\P^{1, an}_{v}\simeq \P^1(\C)$. Thus, we may work on $\P^1(\C)$. For each $n\in\N$ we let 
  $$G_{F, C, n, v}(t_1, t_2):=\frac{\log \| F_{C, n}(t_1, t_2)\|_v}{\deg F_{C, n}},$$
which is a plurisubharmonic function on $\C^2\backslash \{(0, 0)\}$. Denote by $\pi$ is the standard projection from $\C^2\backslash \{(0,0)\}$ to $\P^1(\C)$.  As $G_{F, C, n, v}$ converges locally uniformly on $\C^{2}\setminus\{ (t_1, t_2): t_1/t_2\in \Sing(f) \textup{ or } t_1{=}t_2{=}0\}$; see \cite{Branner:Hubbard:1}, the sequence  $\mu_{n, v}:=\pi_*dd^cG_{F, C, n, v}$ of probability measures 
converges weakly  to the rescaled bifurcation measure $\mu_c/\hhat_f(c)$ on $\P^1(\C)\backslash \Sing(f)$. Since the space of probability measures on $\P^1(\C)$ is compact in the weak topology, to show that $\mu_{n, v}$ has a unique limit, it suffices to prove that for any convergent subsequence of $\mu_{n, v}$, the limit admits no point mass on $\Sing(f)$. Without loss of generality, we may assume that $0\in \Sing(f)$. We have to show that for any $\epsilon>0$, there is a radius $r>0$ and an integer $N>0$, such that for all $n\geq N$ we have 
  $\mu_{n, v}(D(0, r))<\epsilon.$
Suppose that this is not the case. Then we may find integers $n_j\to \infty$ and a sequence of radii $r_{n_j}\to 0$ such that
   $$\mu_{n_j, v}(D(0, r_{n_j}))\to \epsilon_0>0,$$
 as $j\to\infty$. Let $P_{n_j}(t)$ be a potential function of $\mu_{n_j, v}|_{D(0, r_{n_j})}$. We have
  $$P_{n_j}(t):=\int \log |t-s|_vd(\mu_{n_j, v}|_{D(0, r_{n_j})})\to \epsilon_0\log|t|_v$$
locally uniformly on a punctured disk centered at $0$. Hence the sequence of subharmonic functions $G_{F, C, n_j, v}(t, 1)-P_{n_j}(t)$ converges locally uniformly to a subharmonic function $G_{F, C, v}(t, 1)-\epsilon_0 \log |t|_v$ on a punctured disk. So we can find some $L_0>0$ and $r_0>0$, such that for all big  $n_j$ we have
    $$\sup_{|t|=r_0}\left(G_{F, C, n_j, v}(t, 1)-P_{n_j}(t)\right)<L_0.$$
From \cite[Proposition 3.1]{DeMarco: heights}, one has
  $G_{F, C, v}(t, 1)=o(\log|t|_v).$
Then for very small $t$, we can find $n_j$ big enough such that 
  $$G_{F, C, n_j, v}(t, 1)-P_{n_j}(t)>-\frac{\epsilon_0}{2} \log|t|_v>L_0,$$
which is a contradiction as the subharmonic function $G_{F, C, n_j, v}(t, 1)-P_{n_j}(t)$ achieves its maximal value on the boundary of $D(0, r_0)$. Since all $G_{F, C, n, v}$ are bounded above uniformly near $0$,  $G_{F, C,v}(t, 1)$ is bounded above and subharmonic on the punctured disk centered at $0$, by \cite[Theorem 3.6.1]{R:potential}, $G_{F, C, v}(t,1)$ has a unique extension to a subharmonic function in a disk centered at $0$, with $G_{F, C, v}(0,1):=\lim\sup_{t\to 0} G_{F, C, v}(t,1)$. Because $G_{F, C, v}(t, 1)=o(\log|t|_v)$ and $\mu_{f, c, v}(0)=0$, the extended subharmonic function is a potential of $\mu_{f, c, v}$ near $0$. Hence $\mu_{f, c, v}$ has continuous potential if and only if $G_{F, C, v}(t_1, t_2)$ can be extended continuously. For properties of subharmonic functions, we refer the reader to the book \cite{R:potential}.

Assume now that $v\cM_k$ is non-archimedean. Each $F_{C, n, v}$ determines a probability measure $\mu_{n, v}$ with continuous potential on $\P^{1, an}_{v}$ defined as 
  $$g_{n, v}(x) =\log\|\tilde{x}\|_v- G_{F, C, n, v}(\tilde{x}),$$
which is well defined on $\P^1(\C_v)$ and extends continuously to $\P^{1, an}_{v}$ with
\begin{equation}\label{measurepotential}\mu_{n, v}:=\triangle g_{n, v}+\lambda_v.
\end{equation}
Here $\lambda_v$ is the probability measure supported on the Gauss point. For any neighborhood $U\subset \P^1(\C_v)$ of $\Sing(f)$, the sequence $G_{F, C, n, v}(\tilde x)$ converges uniformly for $x\in  \P^1(\C_v)\backslash U$. Since $\P^1(\C_v)$ is dense in $\P^{1, an}_{v}$, we also have that for any neighborhood of $U^{an}$ of $\Sing(f)$ in $\P^{1, an}_{v}$, the function $g_{n,v}(x)$ converges uniformly on  $\P^{1, an}_{v}\backslash U^{an}$ as $n\to \infty$. Hence, from \eqref{measurepotential} we see that the limit
   $$g_v(x):=\lim_{n\to \infty}g_{n, v}(x)=\log\|\tilde x\|_v-G_{F, C, v}(\tilde x)$$
is an element of BVD$(\P^{1, an}_{v})$ (see \cite[Definition 5.11]{BRbook}), with 
   $$\mu_{f, c, v}-\lambda_v:=\triangle g_v.$$
It is clear that the probability measure $\mu_{f, c, v}$ is the unique limit of $\{\mu_{n, v}\}$ on $\P^{1, an}_{v}$, with potential $g_v$. Since the potential function of $\mu_{f, c, v}$ is unique up to a constant, we have that $\mu_{f, c, v}$ has a continuous potential if and only if $g_v$ can be extended continuously to $\Sing(f)$, or equivalently if and only if $G_{F, C, v}$ can be extended continuously on $\C_v^2\backslash \{(0, 0)\}$. \qed
\begin{cor}
For a non-isotrivial and non-preperiodic dynamical pair $(f, c)$, the total mass of the bifurcation measure is
   $$\mu_c(\P^1(\C)\backslash\Sing(f))=\hhat_f(c).$$
\end{cor}
We are indebted to Laura DeMarco for sharing the idea of the following proposition. 
\begin{prop}\label{capacity resultant}
Let $(f,c)\in K(z)\times K$ be a non-isotrivial and non-preperiodic dynamical pair and $v\in\cM_k$. We write $M_{F,C,v}=\{(t_1, t_2)\in \C_v^2\setminus\{(0,0)\}~:~ G_{F,C,v}(t_1, t_2)\leq 0\}$. Suppose that $G_{F,C,v}$ extends continuously on $\C_v^2\backslash\{(0, 0)\}$. Then, 
   $$ \capacity(M_{F,C,v})\leq \liminf_{n\to \infty} |\Res F_{C, n}|_v^{-1/\deg(F_{C,n})^2}.$$
\end{prop}
\proof Let $v\in\cM_k$ and write $G_{F, C, n, v}(t_1, t_2):=\frac{\log \|F_{C,n}(t_1, t_2)\|_v}{\deg F_{C,n}}$, and   
  $$M_{n,v}:=\left\{(t_1, t_2)\in \C_v^2\setminus\{(0,0)\}~:~G_{F, C, n, v}(t_1, t_2) \leq 0\right\}.$$
From \cite[Proposition 2.1]{DWY:Lattes}, we have $\capacity(M_{n,v})=|\Res(F_{C,n})|_v^{-1/\deg (F_{C,n})^2}$. We are going to prove that for any $\epsilon\in |\C_v^{*}|_v$, there exists an $N\in\N$ such that for any $n\geq N$ we have
\begin{equation}\label{upper limit}  G_{F, C, n, v}(t_1, t_2)-G_{F,C,v}(t_1,t_2) < \epsilon.
\end{equation}
This will imply that $M_{F,C,v}\subset e^{\epsilon} M_n$ and hence by the monotonicity of the homogeneous capacity $\capacity(M_{F,C,v})\leq  e^{2\epsilon} \capacity(M_{n,v})$ for all $n\ge N$. 
Since $|\C_v^{*}|_v$ is dense in $\R_{\ge 0}$, the proposition follows. Note that $G_{F, C, n, v}(t_1, t_2)$ converges locally uniformly to $G_{F,C,v}(t_1,t_2)$ away from $\Sing(f)$. Hence, it suffices to prove that \eqref{upper limit} holds in a small neighborhood of $\Sing(f)$. To this end, we may assume without loss of generality that $0\in\Sing(f)$ and show that there exists an $r>0$, such that for all $t\in\C_v$ with $|t|_v\leq r$ we have $G_{F,C,n,v}(t,1)< G_{F,C,v}(t,1)+\epsilon$ for large $n$.
Since $G_{F,C,v}(t, 1)$ is continuous, we can choose $r$ small enough such that for $|t|_v\leq r$ we have $|G_{F,C,v}(t, 1)-G_{F,C,v}(0, 1)|<\epsilon/3$. 
Moreover, enlarging $N$ if necessary, we may further assume that $G_{F,C,n,v}(t,1)<G_{F,C,v}(t,1)+\epsilon/3<G_{F,C,v}(0,1)+2\epsilon/3$, when $|t|_v=r$. Then, since $G_{F, C, n, v}(t, 1)$ is subharmonic, by the maximum principle (see \cite[Proposition 8.14]{BRbook} when $v$ is non-archimedean), we get
   $$G_{F, C, n, v}(t, 1)< G_{F,C,v}(0,1)+2\epsilon/3<G_{F,C,v}(t,1)+\epsilon$$
for all $t\in\C_v$ with $|t|_v\leq r$ and $n\ge N$, as claimed. \qed 
   
\begin{definition} We call a non-preperiodic and non-isotrivial dynamical pair $(f,c)$ {\em adelic} or {\em quasi-adelic} if  the corresponding measure $\mu_{f,c}=\{\mu_{f, c, v}\}_{v\in \cM_k}$ defined in Proposition \ref{measure for pair}, is adelic or quasi-adelic respectively. 
\end{definition}

\subsection{A generic dynamical pair is not adelic} 
Let $k$ be a number field or the function field of a smooth projective curve defined over a field of characteristic zero and let $\alpha\in \kbar$ (or $\alpha=\infty \in \P^1$).
In what follows we write
\begin{align}\label{defstar}
\begin{split}
A_{C,n,\al}^{*}(t_1,t_2)&:=\frac{A_{C,n}(t_1,t_2)}{u_{\al}(t_1,t_2)^{m_{C,n}(\al)}}\\
B_{C,n,\al}^{*}(t_1,t_2)&:=\frac{B_{C,n}(t_1,t_2)}{u_{\al}(t_1,t_2)^{m_{C,n}(\al)}}.
\end{split}
\end{align}
If $\alpha\in \kbar$ we use $P_\alpha(z, w), Q_\alpha(z,w)$ to denote $P_{\alpha, 1}(z, w)$ and $Q_{\al, 1}(z,w)$ respectively. In this notation $P_{\infty}(z,w):=P_{1,0}(z,w)$ and $Q_{\infty}(z,w):=Q_{1,0}(z,w)$. We define 
  $$R_{\al}(z,w):=\gcd (P_{\al}(z, w),Q_{\al}(z,w)),$$
with $R_{\al}(z,1)$ being a monic polynomial.
Furthermore, we write
$$P_{\al}=R_{\al}\cdot P_{\al}^{*}\text{ and }Q_{\alpha}=R_{\al}\cdot Q_{\alpha}^{*}\text{, where }\gcd(P_{\al}^{*},Q_{\al}^{*})=1.$$
Note that $R_{\al}\neq 1$ if and only if $\al\in\Sing(f)$. 
We write 
\begin{align*}
\mathcal{Z}(R_{\al}):=\{t_1/t_2\in \P^1(\overline{k})~:~R_{\al}(t_1,t_2)=0\}.
\end{align*}
Next we define our notion of an $\al$-\emph{generic} dynamical pair $(f,c)$ and subsequently state our theorem. 
\begin{definition}\label{generic}
Let $(f,c)\in k(t)(z)\times k(t)$ be a dynamical pair with $d=\deg_zf\ge 3$ and $\al\in\Sing(f)$. We say that $(f,c)$ is \emph{$\al$-generic} if the following properties are satisfied:
\begin{enumerate}[label=\emph{(P\arabic*)},itemindent=*]
\item \label{pdegatleast2}
$\deg(f_{\alpha})\ge 2$.
\item \label{pnotpoly}
There exists $\rho\in\mathcal{Z}(R_{\al})$ that is not a totally ramified fixed point of $f_{\al}^2$.
\item \label{pgoodgrowth}
For all $n\in\N$ we have $f^n_{\al}(c(\al))\notin\mathcal{Z}(R_{\al})$.
\item \label{pnotprep}
$\hat{h}_{f_{\al}}(c(\al))\neq 0$. 
\end{enumerate}
\end{definition}
Recall that when $k$ is a number field condition \emph{\ref{pnotprep}} is equivalent with $c(\al)$ not being preperiodic for $f_{\al}$. If on the other hand $k$ is a function field, it is also the case that $\hat{h}_{f_{\al}}(c(\al))= 0$ when the pair $(f_{\al},c(\al))$ is isotrivial. 

\begin{theorem}\label{genericnoadelic}
Let $k$ be a number field or the function field of a smooth projective curve defined over a field of characteristic zero. Consider $f\in k(t)(z)$ and $c\in k(t)$  such that $\deg f\ge 3$ and the dynamical pair $(f,c)$ is non-preperiodic and non-isotrivial. If there is an $\al\in\Sing(f)$ such that $(f,c)$ is $\al$-generic, then $(f,c)$ is not adelic.
\end{theorem}
Before we proceed to the proof, we need some preliminary results. 
Let $S\subset \cM_k$ be a finite set containing all the archimedean places. We denote the set of $S$-integers  of $k$ by
$$\mathcal{O}_{S,k}:=\{\alpha\in k~:~|\alpha|_v\le 1\text{ for all }v\notin S\}.$$
If $\phi\in k(z)$ is a rational map and $\al\in\kbar$, we denote the orbit of $\al$ under the action of $\phi$ as
$$\mathcal{O}_{\phi}(\al)=\{\phi^n(\al)~:~n\in\N\}.$$
The following theorem will play a crucial role in our proofs. We thank Patrick Ingram for referring us to it.
\begin{theorem}\cite[Theorem 2.2]{Silverman: roth}\label{roth}
Let $\phi\in k(z)$ be a rational map of degree at least $2$ such that  $\phi^2(z)\notin k[z]$ and let $\al\in k$. Let $S\subset \cM_k$ be a finite set containing all the archimedean places. Then $|\mathcal{O}_{\phi}(\al)\cap \mathcal{O}_{S,k}|<\infty.$
\end{theorem}
We point out here that an analog of Theorem \ref{roth} also holds for function fields of curves over a field of characteristic zero; see \cite[Theorem 1]{intff}. 
\begin{lemma}\label{manyprimes}
Let $\phi\in k(z)$ be a rational map of degree at least $2$ such that  $\phi^2(z)\notin k[z]$ and let $(\al,\beta)\in k^2\setminus\{(0,0)\}$ be such that 
$\hat{h}_{\phi}(\frac{\al}{\beta})\neq 0$. Let $\{a_n\}$ and $\{b_n\}$ be the sequences defined as follows. For any choice of coprime homogeneous polynomials $P,Q\in k[z,w]$ such that $\phi=[P:Q]$, we let
\begin{align*}
a_0&=\al \text{ and } b_0=\beta \text{, and for all }n\ge 0
\end{align*}
\begin{align*}
a_{n+1}&=P(a_n,b_n), \text{ }b_{n+1}=Q(a_n,b_n).
\end{align*} 
Then there are infinitely many non-archimedean places $v\in\cM_k$ such that $|b_n|_v<1$ for some $n\in\mathbb{N}$.
\end{lemma}
\begin{proof}
We assume that the statement is false and then derive a contradiction. There exists a finite set $S\subset\cM_k$ containing the archimedean places such that $|b_n|_v\geq 1$ for all $v\notin S$ and all $n\in\mathbb{N}$. We may enlarge the set $S$ if necessary to assume that the coefficients of $P$ and $Q$ are in $\mathcal{O}_{S,k}$ and that for all $v\notin S$ we have $\max\{|\al|_v,|\beta|_v\}=1$. This implies that for all $v\notin S$ and all $n\in\N$ we have $\max\{|a_n|_v,|b_n|_v\}\le 1$. Combining this with our hypothesis we get that $|b_n|_v=1$ and $|a_n|_v\le 1$ for all $v\notin S$ and all $n\in\mathbb{N}$. Therefore $\phi^n(\frac{\al}{\beta})=\frac{a_n}{b_n}\in \mathcal{O}_{S,k}$ for all $n\in\N$. 
Since $\hat{h}_{\phi}(\frac{\al}{\beta})\neq 0$ , in view of Theorem \ref{roth} (see also \cite[Theorem 1]{intff} for the function field case) we get that $\phi^2\in k[z]$.
This contradicts our assumption and concludes the proof.
\end{proof}

\begin{lemma}\label{goodgrowth}
Let $\al\in\P^1(\overline{k})$. Assume that for all $n\in\N$ we have $f_{\al}^n(c(\al))\notin\mathcal{Z}(R_{\al})$. Then $m_{C,n+1}(\al)=d\cdot m_{C,n}(\al)$ for all $n\in\N$. In particular,
\begin{align*}
A_{C,n+1,\al}^{*}(t_1,t_2)&=P_{t_1,t_2}(A_{C,n,\al}^{*}(t_1,t_2),B_{C,n,\al}^{*}(t_1,t_2))\\
B_{C,n+1,\al}^{*}(t_1,t_2)&=Q_{t_1,t_2}(A_{C,n,\al}^{*}(t_1,t_2),B_{C,n,\al}^{*}(t_1,t_2)).
\end{align*}
\end{lemma}

\begin{proof}
Let $\al\in\P^1(\overline{k})$. We will see that $m_{C,n+1}(\al)=d\cdot m_{C,n}(\al).$ In view of \eqref{recursionnot} and \eqref{defstar}, we have
\begin{align*}
\ord_{\al}A_{C,n+1}(t_1, t_2)&=d\cdot m_{C,n}(\al)+\ord_{\al}P_{t_1,t_2}(A_{C,n,\al}^{*}(t_1,t_2),B_{C,n,\al}^{*}(t_1,t_2))\\
\ord_{\al}B_{C,n+1}(t_1, t_2)&=d\cdot m_{C,n}(\al)+\ord_{\al}Q_{t_1,t_2}(A_{C,n,\al}^{*}(t_1,t_2),B_{C,n,\al}^{*}(t_1,t_2)).
\end{align*}
We claim the $u_\al$ does not divide both $P_{t_1,t_2}(A_{C,n,\al}^{*},B_{C,n,\al}^{*})$ and $Q_{t_1,t_2}(A_{C,n,\al}^{*},B_{C,n,\al}^{*})$.
To see this, assume the contrary. Then
$$P_{\al}(A_{C,n,\al}^{*}(\al),B_{C,n,\al}^{*}(\al))=Q_{\al}(A_{C,n,\al}^{*}(\al),B_{C,n,\al}^{*}(\al))=0.$$
This implies that either $R_{\al}(A_{C,n,\al}^{*}(\al),B_{C,n,\al}^{*}(\al))=0$, contradicting our assumption that $f_{\al}^n(c(\al))=\frac{A_{C,n,\al}^{*}(\al)}{B_{C,n,\al}^{*}(\al)}\notin \mathcal{Z}(R_{\al})$, or
   $$P^{*}_{\al}(A_{C,n,\al}^{*}(\al),B_{C,n,\al}^{*}(\al))=Q^{*}_{\al}(A_{C,n,\al}^{*}(\al),B_{C,n,\al}^{*}(\al))=0,$$
contradicting the fact that $P^{*}_{\al}$ and $Q^{*}_{\al}$ have no common factor. Here we note that \eqref{defstar} implies that $(A_{C,n,\al}^{*}(\al),B_{C,n,\al}^{*}(\al))\neq (0,0)$. Therefore $m_{C,n+1}(\al)=d\cdot m_{C,n}(\al),$ as claimed. 
The rest of the lemma now follows by \eqref{recursionnot}.
\end{proof}

To state the next proposition, leading to our main theorem of this section, we mention that a rational map $\phi:\P^1\to\P^1$, defined over an algebraic closed field $k$ of characteristic zero, is called a \emph{polynomial map} if it has a totally ramified fixed point, that is, if there is $\al\in\P^1$ such that $\phi^{-1}(\{\alpha\})=\{\alpha\}$. It is well known that if $\phi$ is of degree $d\ge2$ and $\phi^n$ is a polynomial map for some $n\geq 2$, then already $\phi^2$ is a polynomial map. Moreover, if $\phi$ is not a polynomial but $\phi^2$ is a polynomial, then $\phi$ is linearly conjugate to $\frac{1}{z^d}$; see \cite[Proposition 1.1]{Silverman: roth}.

\begin{proposition}\label{theprimes}
Let $(f,c)\in k(t)(z)\times k(t)$ be a dynamical pair with $\deg f\ge 3$ and $\alpha \in\Sing(f)$ be such that $(f,c)$ is $\al-$generic. Then there are infinitely many $v\in\cM_{k}$ such that for some $n_v\in\N$ we have 
$$\max\{|A_{C,n_v,\al}^{*}(\al)|_v,|B_{C,n_v,\al}^{*}(\al)|_v\}<1. $$
\end{proposition}
\begin{proof} 
Let $\alpha \in\Sing(f)$ be such that $(f,c)$ is $\al-$generic. To simplify the notation, throughout this proof we write $a_n^{*}:=A_{C,n,\al}^{*}(\al)$ and $b_n^{*}:=B_{C,n,\al}^{*}(\al)$.
Since \emph{\ref{pgoodgrowth}} holds, Lemma \ref{goodgrowth} yields
\begin{align}\label{recurseventually}
\begin{split}
a^{*}_{n+1}&=P_{\al}(a^{*}_n,b^{*}_n)=R_{\al}(a^{*}_n,b^{*}_n)P^{*}_{\al}(a^{*}_n,b^{*}_n)\\
b^{*}_{n+1}&=Q_{\al}(a^{*}_n,b^{*}_n)=R_{\al}(a^{*}_n,b^{*}_n)Q^{*}_{\al}(a^{*}_n,b^{*}_n),
\end{split}
\end{align}
for all $n\in\N$. We define auxiliary sequences $\{a_n\}$ and $\{b_n\}$ as $a_0=a^{*}_0, b_0=b^{*}_0$ and 
\begin{align*}
   a_{n+1}=P^{*}_{\al}(a_n,b_n), ~ b_{n+1}=Q^{*}_{\al}(a_n,b_n)\text{ for } n\ge 1.
\end{align*}
Then, 
 for all $n\in\N$ we have
\begin{align}\label{enoughforanbn}
\begin{split}
a^{*}_{n}=\displaystyle\prod_{i=0}^{n-1}R_{\al}(a_i,b_i)^{d^{n-1-i}}\cdot a_n\\
b^{*}_{n}=\displaystyle\prod_{i=0}^{n-1}R_{\al}(a_i,b_i)^{d^{n-1-i}}\cdot b_n.
\end{split}
\end{align}
Let $S\subset\cM_{k}$ be a finite set of places, containing the archimedean ones, such that for all $v\notin S$,  the coefficients of $P_{\al}^{*}$ and $Q_{\al}^{*}$ are $v$-adic integers, $|\mathrm{Res}(P_{\al}^{*},Q_{\al}^{*})|_v=1$ and $\max\{|b_0|_v,|a_0|_v\}=1$. 
Then invoking \cite[Lemma 10.1]{BRbook}, for all $n\in\mathbb{N}$ and $v\notin S$ we have 
\begin{align}\label{Coprime}
\max\{|a_n|_v,|b_n|_v\}=1.
\end{align}
We may enlarge the set $S$ if necessary to assume that the elements of $\mathcal{Z}(R_{\al})\cap \kbar$ are $v$-adic integers for all $v\notin S$. Now
combining \eqref{enoughforanbn} with \eqref{Coprime} we get that 
\begin{align}\label{reducetochangeofvariables}
\max\{|a^{*}_{n+1}|_v,|b^{*}_{n+1}|_v\}\le |u_{\rho}(a_n,b_n)|_v,
\end{align}
for all $v\notin S$ and $n\in\N$ and for any $\rho\in\mathcal{Z}(R_{\al})$.
Now let $\rho\in\mathcal{Z}(R_{\al})$ be as in \emph{\ref{pnotpoly}}. We claim that there are infinitely many $v\in\cM_{k}$ such that
\begin{align}\label{primitive}
|u_{\rho}(a_{n}, b_{n})|_v<1,
\end{align}
for some $n\in\mathbb{N}$. By \eqref{reducetochangeofvariables} it is clear that this suffices to prove this proposition.
To prove \eqref{primitive}, we use Lemma \ref{manyprimes}. If $\rho=\infty$ our claim follows. Otherwise, let $M_{\rho}(z,w)=(w+\rho z,z)$ and $(\hat{P_{\al}^{*}},\hat{Q_{\al}^{*}})=M_{\rho}^{-1}\circ(P_{\al}^{*},Q_{\al}^{*})\circ M_{\rho}$. Consider the morphism $\hat g_{\al}:\P^1\to\P^1$ defined by
\begin{align*}
[z:w]\mapsto[\hat{P}_{\al}^{*}(z,w):\hat{Q}_{\al}^{*}(z,w)].
\end{align*}
Since $\rho\in\mathcal{Z}(R_{\al})$ is as in \emph{\ref{pnotpoly}}, we know that $\infty$ is not a totally ramified fixed point of $\hat g^2_{\al}$.
Moreover, by \emph{\ref{pnotprep}} we have $|\mathcal{O}_{\hat g_{\al}}(\frac{b_0}{a_0-\rho b_0})|=\infty$ and by  \emph{\ref{pdegatleast2}} we have $\deg(g_{\al})\ge 2$.
Thus, by Lemma \ref{manyprimes} applied to the rational map $\hat g_{\al}$ and  $(b_0,a_0-\rho b_0)$ our claim follows. This finishes our proof.
\end{proof}
We are now ready to prove the main result of this section: In most cases $(f,c)$ is not adelic.
\subsection{Proof of Theorem \ref{genericnoadelic}} 
 Let $\alpha\in\Sing(f)$ be such that $(f,c)$ is $\al-$generic and assume to the contrary that $(f,c)$ is adelic. Then there is a finite set $S\subset\cM_k$ such that for all $v\notin S$ we have
\begin{align}\label{adelic}
G_{F, C, v}(t_1,t_2)=\log \|(t_1,t_2)\|_v+c_v,
\end{align}
for a constant $c_v$. 
Denote by $\mathcal{P}\subset\cM_k$ the infinite set of places satisfying the conclusion of Proposition \ref{theprimes}. In other words, for $v\in\mathcal{P}$, there exists $n_v\in\N$ such that
$$\max\{|A_{C,n_v,\al}^{*}(\al)|_v,|B_{C,n_v,\al}^{*}(\al)|_v\}<1.$$
Enlarging the set $S$ if necessary we may further assume that for all $v\notin S$, the following hold
\begin{enumerate}[label=\emph{(S\arabic*)},itemindent=*]
\item\label{aisunit}
$|\al|_v=1$ if $\al\neq \infty$.
\item \label{forgcdabs1}
For all $\beta\in\Sing(f)\setminus\{\al\}$ we have $|u_\beta(\al)|_v=|u_\beta(0,1)|_v=1$.
\item\label{vintegrallift}
The coefficients of $F_{t_1,t_2}$ are $v-$adic integers.
\item \label{forgoodreduction}
$|\Res(A,B)|_v=1$ and the coefficients of $C(t_1,t_2)=(A, B)$ are $v-$adic integers.
\item \label{reshasunitcoef}
All coefficients of $\Res_{(z,w)}(F_{t_1,t_2})$ are $v$-adic units.
\end{enumerate}
We aim to prove that \eqref{adelic} does not hold for places in the infinite set $\mathcal{P}\setminus S$, thus leading to a contradiction. To do so we will evaluate \eqref{adelic} at two distinct points that yield distinct values for $c_v$ when $v\in\mathcal{P}\setminus S$.

Let $v\in\mathcal{P}\setminus S$. In the rest of this proof, for $t_0\in \C_v$ we write
\[ T_0 = \left\{
\begin{array}{ll}
(t_0,1) &\mbox{ if $\al\in\A^1(\overline{k})$}\\
(1,t_0) &\mbox{ if $\al=\infty$.}\\
\end{array} \right. \]
View both $T_0$ and $\alpha$ as elements of $\P^1$. Since either $A_{C,n_v,\al}^{*}(\al)$ or $B_{C,n_v,\al}^{*}(\al)$ is non-zero and $v$ is non-archimedean,  we may choose $T_0$ sufficiently close to $\alpha$ in the $v-$adic topology to be such that
\begin{enumerate}[label=\emph{(T\arabic*)},itemindent=*]
\item \label{veryclose}
$0<|u_{\al}(T_0)|_v<1$.
\item \label{fromttoa}
$\max\{|A_{C,n_v,\al}^{*}(T_0)|_v,|B_{C,n_v,\al}^{*}(T_0)|_v\}\le \max\{|A_{C,n_v,\al}^{*}(\al)|_v,|B_{C,n_v,\al}^{*}(\al)|_v\}$.
\end{enumerate}
Next, we show that evaluating \eqref{adelic} at $T_0$ gives $c_v<0$. To this end, recall that 
\begin{align}\label{overgcdoftherest}
F_{C,n}(t_1,t_2)=\left(\frac{A^{*}_{C,n,\al}(t_1,t_2)}{g_{C,n,\al}(t_1,t_2)},\frac{B^{*}_{C,n,\al}(t_1,t_2)}{g_{C,n,\al}(t_1,t_2)}\right),
\end{align}
where
$$g_{C,n,\al}(t_1,t_2)=\prod_{\substack{\beta\in\Sing(f)\setminus\{\al\}}}u_{\beta}(t_1,t_2)^{m_{C,n}(\beta)}.$$
Combining  \emph{\ref{forgcdabs1}} and \emph{\ref{veryclose}}, the ultrametric inequality gives $|u_\beta(T_0)|_v=1$ for all $\beta\in\Sing(f)\setminus\{\al\}$. This in turn yields $|g_{C,n,\al}(T_0)|_v=1$ for all $n\in\N$. Thus, evaluating \eqref{overgcdoftherest} at $T_0$, we have
\begin{align*} 
\|F_{C,n}(T_0)\|_v=\|\left(A^{*}_{C,n,\al}(T_0),B^{*}_{C,n,\al}(T_0)\right)\|_v.
\end{align*}
Moreover, since \emph{\ref{pgoodgrowth}} holds, we can use Lemma \ref{goodgrowth} to get that for all $n\ge n_v$ we have
\begin{align}\label{primegrowth}
\|F_{C,n+1}(T_0)\|_v=\|F_{T_0}(F_{C,n}(T_0))\|_v.
\end{align}
Note that by \emph{\ref{aisunit}} and \emph{\ref{veryclose}} we have $||T_0||_v\le 1$. Combining this with \emph{\ref{vintegrallift}} we get that for all $n\ge n_v$ the following inequality follows from \eqref{primegrowth}:
$$\|F_{C,n+1}(T_0)\|_v\le \|F_{C,n}(T_0)\|_v^{d}.$$
An easy argument by induction and \emph{\ref{fromttoa}} yield that for for all $n\ge n_v$ we have
$$\|F_{C,n}(T_0)\|_v\le\max\{|A_{C,n_v,\al}^{*}(\al)|_v,|B_{C,n_v,\al}^{*}(\al)|_v\}^{d^{n-n_v}}<1,$$
where the last inequality follows from our assumption that $v\in\mathcal{P}$.
We now get
\begin{align}\label{resultantdivisible}
c_v=\lim_{n\to\infty}\frac{\log\|F_{C,n}(T_0)\|_v}{\deg(F_{C,n})}\le\frac{\log\max\{|A_{C,n_v,\al}^{*}(\al)|_v,|B_{C,n_v,\al}^{*}(\al)|_v\}}{d^{n_v}\cdot \hat{h}_f(c)}<0,
\end{align}
as claimed. We point out here that by our assumption the dynamical pair $(f,c)$ is not isotrivial. Hence, our property \emph{\ref{pnotprep}} guarantees that  $\hat{h}_f(c)\neq 0$; see \cite{Baker, DeMarco: heights}. 

On the other hand, we can choose $S_0=(s_0,1)\in\C^2_v$ to be such that for all $\beta\in\Sing(f)$ we have $|s_0|_v=|u_\beta(S_0)|_v=1$. Then, upon using \emph{\ref{vintegrallift}}, we get that the coefficients of $F_{S_0}$ are $v-$adic integers. Moreover, by \emph{\ref{reshasunitcoef}} and since $|u_\beta(S_0)|_v=1$ for all $\beta\in\Sing(f)$, we have $|\Res_{(z,w)}(F_{S_0})|_v=1$. Therefore, \cite[Lemma 10.1]{BRbook} yields that $\|F_{S_0}(z,w)\|_v=\|(z,w)\|_v^d$.
Since by \emph{\ref{forgoodreduction}} we have that $C$ has good reduction and moreover $\|S_0\|_v=1$, another application of \cite[Lemma 10.1]{BRbook} yields $\|C(S_0)\|_v=1$. Now inductively we have $\|F^n_{S_0}(C(S_0))\|_v=1$. By our choice of $S_0$ we have $|g_{C,n}(S_0)|_v=1$ for all $n\in\N$. Hence $\|F^n_{S_0}(C(S_0))\|_v=\|F_{C,n}(S_0)\|_v=1$ for all $n\in\N$. Therefore,
\begin{align*}
c_v=\lim_{n\to\infty}\frac{\log\|F_{C,n}(S_0)\|_v}{\deg(F_{C,n})}=0.
\end{align*}
This contradicts \eqref{resultantdivisible} our assumption and finishes the proof of our theorem.
\qed

\begin{remark}\label{moregeneric}
A dynamical pair $(f,c)\in k(t)(z)\times k(t)$ is adelic if and only if the pair $(M\circ f^n\circ M^{-1},M(f^N(c)))$ is adelic for some $n,N\in\N$ and a M\"obius transformation $M(z)\in\overline{k(t)}(z)$. 
\end{remark}

Theorem \ref{cornoadelic} is a direct corollary of Theorem \ref{genericnoadelic} and Remark \ref{moregeneric}, once one notices that conditions \emph{\ref{pgoodgrowth}} and \emph{\ref{pnotprep}} will be satisfied as long as $\hhat_{f_\al}(c(\al))$ is strictly bigger than $L:=\max\{\hhat_{f_\al}(a)~:~a\in\mathcal{Z}(R_{\al})\}$. The later is a well defined quantity since the set  $\mathcal{Z}(R_{\al})$ is finite and is independent of the starting point $c\in k(t)$.
\qed

\section{Quasi-adelicity for almost all starting points }\label{quasi-adelic family}

We study the family $g_{\lambda, t}(z):=\frac{\lambda z}{z^2+tz+1}$, where $\lambda$ is an $\ell$-th primitive root of unity for $\ell\geq 2$, and aim to prove Theorem \ref{mainthm3}. We show that for a generic $c$, the dynamical pair $(g_{\lambda, t}, c)$ is quasi-adelic but is not adelic. In particular, our theorem can be applied with $c\in\{1,-1\}$ being a critical point of $g_{\lambda, t}$. We refer the reader to \cite{Milnor:quad, BEE} for the pictures of the bifurcation of $(g_{\lambda, t}(z), \pm 1)$. 
\subsection{Homogenous lifts} Throughout the rest of this section, we fix a primitive $\ell$-th root of unity $\lambda$ with order at least $2$. It is more convenient to work with the $\ell$-th iterate of $g_{\lambda, t}(z)$, which we denote by
   $$f_t(z):=g_{\lambda, t}^{\ell}(z)$$
Since $g_{\lambda, t}$ degenerates at $t=\infty$, we have $\Sing(f)=\{\infty\}$. Let us now introduce some notation that will be used throughout this section; we write
  $$d:=2^{\ell}, ~ d_1:=2^{\ell-1}-1, ~d_2:=2^{\ell-1}.$$
At times we also use notation introduced in Section \S\ref{a pair is rarely adelic}.
We fix a homogenous lift of $g_{\lambda, t}(z)$ as
   $$G_{t_1,t_2}(z, w):=(t_2\lambda zw, t_1zw+t_2(z^2+w^2)).$$
We begin with establishing a proposition which enables us to show that $g^\ell_{\lambda, t}=f_t$ degenerates to $f_\infty(z):=z/(z^2+1)$ at $t=\infty$. 
\begin{prop}\label{degenerations} Let $\lambda$ be an $\ell$-th primitive root of unity. The $\ell$-th iterate of $G_{t_1,t_2}$ is given by 
  $$G^{\ell}_{t_1,t_2}(z, w)=t_2^{d_2}\cdot \left(\tau \cdot (t_1zw)^{d_1}zw+t_2(\cdots), \tau \cdot (t_1zw)^{d_1} (z^2+w^2)+t_2(\cdots)\right),$$
for some non-zero $\tau\in \Z[\lambda]$. In particular, $f_t=g^\ell_{\lambda, t}$ degenerates to $f_\infty(z)=z/(z^2+1)$ at $t=\infty$. 
 \end{prop}
 \proof We prove this proposition by induction. Notice that, inductively for $2\leq n\leq \ell$ one has $G^n_{t_1, t_2}=(P_n, Q_n)$, where 
\begin{align}\label{Gl}
\begin{split}
P_n&=t_2^{2^{n-1}}\cdot\left((zw)^{2^{n-1}}\alpha_n t_1^{2^{n-1}-1}  +(zw)^{2^{n-1}-1}(z^2+w^2)\eta_{n}t_2t_1^{2^{n-1}-2}+t_2^{2}(\cdots)\right)\\
Q_n&=t_2^{2^{n-1}-1}\cdot\left( (zw)^{2^{n-1}}\beta_{n} t_1^{2^{n-1}} +(zw)^{2^{n-1}-1}(z^2+w^2)\tau_{n} t_2t_1^{2^{n-1}-1}+t_2^{2}(\cdots)\right), 
\end{split} 
\end{align}
for constants $\alpha_{n}, \beta_n, \eta_n$ and $\tau_n$ depending on $\lambda$. From the iteration formula we get $\alpha_2=\lambda^2, \eta_2=\lambda^2, \beta_{2}=1+\lambda,\tau_{2}=\lambda +2$ and for all $n\ge2$ we have
\[ \left\{
\begin{array}{ll}
\alpha_{n+1}=\lambda \cdot \alpha_{n} \cdot \beta_{n},\\
 \beta_{n+1}= \beta_{n}\cdot(\alpha_{n}+ \beta_{n})\\
\end{array} \right. 
\textup{and }\left\{
\begin{array}{ll}
\eta_{n+1}=\lambda \cdot ( \alpha_{n} \cdot \tau_{n}+\beta_n\cdot \eta_n),\\
 \tau_{n+1}= \alpha_n\cdot \tau_n+\beta_n\cdot \eta_n+2\beta_n\cdot \tau_n.\\
\end{array} \right.  \]
Consequently, for $n\geq 3$ we have
\begin{align*} 
\alpha_{n}&=\lambda^n\cdot \prod_{i=1}^{n-2}(1+\lambda+\cdots +\lambda^i)^{2^{n-2-i}},\\
 \beta_{n}&=(1+\lambda+\cdots+\lambda^{n-1})\cdot \prod_{i=1}^{n-2}(1+\lambda+\cdots +\lambda^i)^{2^{n-2-i}}.
\end{align*}
Since $\lambda$ is an $\ell$-th primitive root of unity, we have $\alpha_\ell\neq 0$ and $\beta_\ell=0$. It remains to show that 
\begin{align}\label{ellthequal}   
\tau:=\tau_\ell=\alpha_\ell.
\end{align}
Let $z_0\in \C$ be such that $\frac{z_0}{z_0^2+1}=1$. Since $\beta_\ell=0$,  from the expression of $G^{\ell}_{t_1, t_2}(z_0, 1)$ in \eqref{Gl} we get
   $$\lim_{t_1=1, t_2\to 0}\frac{P_\ell(z_0,1)}{Q_\ell(z_0, 1)}\to \frac{ \alpha_\ell}{\tau_\ell},$$
or equivalently 
\begin{equation}\label{equaltau}
\lim_{t\to \infty} g_{\lambda, t}^\ell(z_0)= \frac{ \alpha_\ell}{\tau_\ell}.
\end{equation}
We are going to show that this limit is equal to one; hence equation \eqref{ellthequal} follows. Notice that 
   $$g_{\lambda, t}(z_0)=\frac{\lambda  z_0}{1+t z_0+z_0^2}=\frac{\lambda}{t}\cdot \left(1-\frac{1}{t}+o\left(\frac{1}{t}\right)\right),$$
for $t\to \infty$. 
Using the expression of $g_{\lambda, t}(z)=\lambda \cdot z/(1+t\cdot z+z^2)$, inductively we get 
   $$g^n_{\lambda, t}(z_0)=\frac{1}{t}\cdot \frac{\lambda^n}{1+\lambda+\cdots +\lambda^{n-1}}\cdot \left(1-\frac{1}{1+\lambda+\cdots+\lambda^{n-1}}\cdot \frac{1}{t}+o\left(\frac{1}{t}\right)\right),$$
for all $1\leq n\leq \ell-1$ as $t\to \infty$. Consequently we have
\[\begin{split}
 g^\ell_{\lambda, t}(z_0)&=g_{\lambda, t}( g^{\ell-1}_{\lambda, t}(z_0))\\
 &=g_{\lambda, t}\left(\frac{1}{t}\cdot \frac{\lambda^{\ell-1}}{1+\lambda+\cdots +\lambda^{\ell-2}}\cdot \left(1-\frac{1}{1+\lambda+\cdots+\lambda^{\ell-2}}\cdot \frac{1}{t}+o\left(\frac{1}{t}\right)\right)\right)\\
 &=\frac{\frac{1}{t}\cdot \frac{\lambda\cdot \lambda^{\ell-1}}{1+\lambda+\cdots +\lambda^{\ell-2}}\cdot \left(1-\frac{1}{1+\lambda+\cdots+\lambda^{\ell-2}}\cdot \frac{1}{t}+o\left(\frac{1}{t}\right)\right)}{1+t\cdot \frac{1}{t}\cdot \frac{\lambda^{\ell-1}}{1+\lambda+\cdots +\lambda^{\ell-2}}\cdot \left(1-\frac{1}{1+\lambda+\cdots+\lambda^{\ell-2}}\cdot \frac{1}{t}+o\left(\frac{1}{t}\right)\right)+o\left(\frac{1}{t}\right)}\\
 &=\frac{\frac{1}{t}\cdot \frac{1}{1+\lambda+\cdots +\lambda^{\ell-2}}\cdot \left(1-\frac{1}{1+\lambda+\cdots+\lambda^{\ell-2}}\cdot \frac{1}{t}+o\left(\frac{1}{t}\right)\right)}{1+(-1)\cdot \left(1-\frac{1}{1+\lambda+\cdots+\lambda^{\ell-2}}\cdot \frac{1}{t}+o\left(\frac{1}{t}\right)\right)+o\left(\frac{1}{t}\right)},
 \end{split}\]
where in the last equality we used the fact that $\lambda^\ell=1$ and $1+\lambda+\cdots +\lambda^{\ell-1}=0$. Letting now $t\to \infty$, we get $g^\ell_{\lambda, t}(z_0)\to 1$. Combining this with \eqref{equaltau} we get \eqref{ellthequal}. The proposition follows.  \qed

Let us now fix a lift of $f_t$ in homogeneous coordinates $(t_1,t_2)$ as
\begin{align}\label{tau}
\begin{split}
F_{t_1,t_2}(z,w):&=(P_{t_1,t_2}(z,w),Q_{t_1,t_2}(z,w)):=G^\ell_{t_1,t_2}(z, w)/(\tau\cdot t_2^{d_2})\\
&= \left( (t_1zw)^{d_1}zw+t_2(\cdots), (t_1zw)^{d_1} (z^2+w^2)+t_2(\cdots)\right).
\end{split}
\end{align}
Notice that at the point at infinity our lift specializes to the map
$$F_{1,0}(z,w)=\left( (zw)^{d_1}zw, (zw)^{d_1} (z^2+w^2)\right),$$
which is a homogenious lift of  
   $$f_{\infty}(z)=\frac{z}{z^2+1}.$$
Keeping the notation as in Section \S\ref{a pair is rarely adelic}, we have $R_\infty(z,w)=(zw)^{d_1}$; hence $\mathcal{Z}(R_{\infty})=\{0,\infty\}$.
Let $k$ be a number field containing $\lambda$, so that $f_t\in k(t)(z)$.
We now fix a starting point $c\in k(t)$ satisfying $\mathbf{0,\infty\notin\mathcal{O}_{f_{\infty}}(c(\infty))}$; compare this with condition \emph{\ref{pgoodgrowth}} in Definition \ref{generic}.
We also fix a homogeneous lift of the starting point $c$, with coefficients in $\mathcal{O}_k$, as
 $$C(t_1,t_2):=(A(t_1,t_2),B(t_1,t_2)),$$
and write
$$F_{C,n}(t_1,t_2)=F^n_{t_1,t_2}(C)/{\rm gcd}(F^n_{t_1, t_2}(C))=(A_{C,n}(t_1,t_2),B_{C,n}(t_1,t_2)).$$

\begin{lemma}\label{gcd}
For all $n\in\N$ we have ${\rm gcd}(F^n_{t_1, t_2}(C))=1$. Hence
$$F_{C,n}(t_1,t_2)=F^n_{t_1,t_2}(C(t_1,t_2)).$$
\end{lemma}

\begin{proof}
As $F_{t_1,t_2}$ only degenerates at $t_1/t_2=\infty$,  we know that ${\rm gcd}(F^n_{t_1, t_2}(C))=t_2^{m_{C,n}(\infty)}$ for all $n\in\N$.
Since $0,\infty\notin\mathcal{O}_{f_{\infty}}(c(\infty))$, Lemma \ref{goodgrowth} yields that $m_{C,n+1}(\infty)=d\cdot m_{C,n}(\infty)$ for all $n\in\N$. Moreover the fact that $c(\infty)\neq 0,\infty$ yields that $m_{C,1}(\infty)=0$. Therefore ${\rm gcd}(F^n_{t_1, t_2}(C))=1$ for all $n\in\N$ and the lemma follows.
\end{proof}

\begin{lemma}\label{degree}
We have $\deg(F_{C,n})=d_1\cdot \frac{d^n-1}{d-1}+d^n\cdot\deg(c)$ for all $n\in\N$. In particular, the dynamical pair $(f, c)$ is not preperiodic. Furthermore  $\hat{h}_f(c)=\frac{d_1}{d-1}+\deg(c)\neq 0$.  
\end{lemma}

\begin{proof}
Since $0,\infty\notin\mathcal{O}_{f_{\infty}}(c(\infty))$ we get  $\deg(A_{C,n})=\deg(B_{C,n})$ for all $n\in\mathbb{N}$. The lemma now follows inductively from the recursive definition of $F_{C,n}$.
\end{proof}

When there is no scope for confusion, in what follows we use $a_n$ and $b_n$ to denote $A_{C,n}(1,0)$ and $B_{C,n}(1,0)$ respectively. From Lemma \ref{gcd}, we see that for all $n\in\N$:
\begin{align}\label{anbnourpair}
\begin{split}
a_{n+1}&=(a_nb_n)^{d_1}\cdot a_nb_n\\
b_{n+1}&=(a_nb_n)^{d_1}\cdot (a^2_n+b^2_n).
\end{split}
\end{align}
We also make use of auxiliary sequences $\{a^*_{n}\}, \{b^*_{n}\}\subset k$ defined by 
$a_0^{*}=a_0$, $b_0^{*}=b_0$ and for $n\ge1$:
\begin{equation}\label{anbnstar}
 a_{n+1}^*=a_n^*b_n^*, ~ b^*_{n+1}={a_n^*}^2+{b_n^*}^2.
\end{equation}
Notice that if we define 
\begin{align}\label{aldef}
\al_n:=\prod_{i=0}^{n-1}(a_i^{*}b_i^{*})^{ d_1\cdot d^{n-i-1}},
\end{align}
then for $n\geq 1$ we have $a_n=\al_na_n^{*}$ and $b_n=\al_nb_n^{*}$. 

\subsection{Continuity of the escape rate}

In order to prove that $(f,c)$ is quasi-adelic, we need to first show that the escape rate $G_{F, C, v}$ is a continuous function. 

\begin{theorem}\label{continuity}
The functions $\frac{\log\|F_{C,n}(t_1,t_2)\|_v}{\deg(F_{C,n})}$ converge locally uniformly on $\C_v^2\setminus\{(0,0)\}$ to the function $G_{F,C,v}$. In particular $G_{F,C,v}$ is continuous.
\end{theorem}

Before we proceed to the proof of this theorem, we establish some lemmata. First we let 
$$F_{C,n}(1,s)=(A_{C,n}(1,s),B_{C,n}(1,s))=(A_{C,n}(1,0)+sp_n(s),B_{C,n}(1,0)+s q_n(s)).$$
\begin{lemma}\label{leadingtermslimit}
For each $v\in\cM_k$, we have the following
\begin{itemize}
\item
$\gamma_v:=\displaystyle\lim_{n\to\infty}\frac{\log|A_{C,n}(1,0)|_v}{d^n}=\lim_{n\to\infty}\frac{\log|B_{C,n}(1,0)|_v}{d^n}$; and 
\medskip
\item
$\displaystyle\limsup_{n\to\infty}\frac{\log|p_n(0)|_v}{d^n}, \displaystyle\limsup_{n\to\infty}\frac{\log|q_n(0)|_v}{d^n}\le \gamma_v.$
\end{itemize}
\end{lemma}

\begin{proof}
By \cite[Lemma 10.1]{BRbook}, the recursive definition of $\{a_n^{*}\},\{b_n^{*}\}$ in \eqref{anbnstar} implies that there is a set of constants $\{ L_v ~:~v\in\cM_k\}$ and a finite set $S\subset\cM_k$ such that $L_v=1\text{ for all }v\notin S$ and for all $v\in\cM_k$ we have $L_v\ge 1$ and
\begin{align}\label{upperboundM}
\max\{|a_n^{*}|_v,|b_n^{*}|_v\}&\le L_v^{2^n}\text{ for all }n\in\N.
\end{align}
Now let $L=\displaystyle\prod_{\substack{v\in\cM_k}}L_v^{N_v}$. Invoking the product formula, inequality \eqref{upperboundM} yields
\begin{align}\label{lowerM}
\min\{|a_n^{*}|_v,|b^{*}_n|_v\}\ge \frac{1}{L^{2^n}}.
\end{align}
Moreover \eqref{upperboundM} implies
\begin{align}\label{upperM}
\max\{|a_n^{*}|_v,|b^{*}_n|_v\}\le L^{2^n}.
\end{align}
In particular \eqref{lowerM} and \eqref{upperM} yield that $\displaystyle\lim_{n\to\infty}\frac{\log|a_n^{*}|_v}{d^n}=\lim_{n\to\infty}\frac{\log|b_n^{*}|_v}{d^n}=0$ and for $\{\al_n\}$ as in \eqref{aldef}, the sequence 
$$\frac{\log|\al_n|_v}{d^n}=\sum_{i=0}^{n-1}d_1\cdot\frac{ \log|a_i^{*}b_i^{*}|_v}{d^{i+1}}$$
converges. Denoting its limit by $\gamma_v$, we have established the following.
$$\gamma_v=\displaystyle\lim_{n\to\infty}\frac{\log|a_n|_v}{d^n}=\lim_{n\to\infty}\frac{\log|b_n|_v}{d^n}=\lim_{n\to\infty}\frac{\log|\al_n|_v}{d^n}.$$ 
The first part of the lemma follows. 
Now let $c_n$ and $e_n$ be the constant terms of the polynomials $p_n(s)$ and $q_n(s)$ respectively. From the recursive definition of $F_{C, n}$ as in Lemma \ref{gcd}, we see that there are homogenous $\Phi, \Psi\in k[z,w]$ of degree $d$ and $\Phi_i, \Psi_i\in k[z,w]$ for $i=1, 2$ of degree $d-1$, such that 
 \begin{align}\label{recursioncd}
\begin{split}
c_{n+1}=\Phi(a_n, b_n)+c_n\cdot \Phi_1(a_n, b_n)+e_n\cdot \Phi_2(a_n, b_n)\\
 e_{n+1}=\Psi(a_n, b_n)+c_n\cdot \Psi_1(a_n, b_n)+e_n\cdot \Psi_2(a_n, b_n),
\end{split}
\end{align}  
for all $n\geq 0$. 
Now we define the sequences $c^{*}_n$ and $e^{*}_n$ as $c_n=\al_n\cdot c^{*}_n$ and $e_n=\al_n \cdot d^{*}_n$.  We are going to show that
\begin{align}\label{wantedstarcd}
\displaystyle\limsup_{n\to\infty}\frac{\log|c^{*}_n|_v}{d^n}\text{, }\displaystyle\limsup_{n\to\infty}\frac{\log|e_n^{*}|_v}{d^n}\le 0.
\end{align}
Having proved this the second part of our lemma will follow, since
\begin{align*}
\displaystyle\limsup_{n\to\infty}\frac{\log|c_n|_v}{d^n}\text{, }\displaystyle\limsup_{n\to\infty}\frac{\log|e_n|_v}{d^n}\le \displaystyle\lim_{n\to\infty}\frac{|\al_n|_v}{d^n}=\gamma_v.
\end{align*}
To prove \eqref{wantedstarcd}, first notice that by \eqref{aldef} we have $\frac{\al^d_n}{\al_{n+1}}=\frac{1}{(a_n^{*}b_n^{*})^{d_1}}$. The recursive formulas in \eqref{recursioncd} can now be written as
\begin{align}\label{reccdstar}
\begin{split}
c^*_{n+1}=\frac{\Phi(a^*_n, b^*_n)+c^*_n\cdot \Phi_1(a^*_n, b^*_n)+e^*_n\cdot \Phi_2(a^*_n, b^*_n)}{(a^*_n b^*_n)^{d_1}}\\
 e^*_{n+1}=\frac{\Psi(a^*_n, b^*_n)+c^*_n\cdot \Psi_1(a^*_n, b^*_n)+e^*_n\cdot \Psi_2(a^*_n, b^*_n)}{(a^*_n b^*_n)^{d_1}}.
\end{split}
\end{align}
Let $L_{n,v}:=\max\{|c^*_n|_v,|e^*_n|_v\}$.
By \eqref{upperboundM} and \eqref{lowerM} we get that there is some constant $L_0\ge 1$ such that
\begin{align*}
\max\left\{|\Phi(a^*_n, b^*_n)|_v,|\Phi_i(a^*_n, b^*_n)|_v,  |\Psi(a^*_n, b^*_n)|_v, |\Psi_i(a^*_n, b^*_n)|_v, \frac{1}{ |a^*_n b^*_n|_v^{d_1}}\right\}\le L_0^{2^n}
\end{align*}
for $i=1, 2$. 
Combining this with \eqref{reccdstar} we see that there is some $r\ge 1$ such that that for all $n\in\N$:
\begin{align*}
L_{n+1,v}\le r^{2^n}L_{n, v}.
\end{align*}
An easy argument by induction now yields that $L_{n,v}\le r^{n2^n}L_{0, v}$
 for all $n\in\N$. Therefore,
\begin{align*}
\displaystyle\limsup_{n\to\infty}\frac{\log L_{n,v}}{d^n}\le 0,
\end{align*}
and  inequalities \eqref{wantedstarcd} follow. This finishes our proof.
\end{proof}

The next two propositions show that the convergence of our escape rate function is locally uniform near the degenerate point at $t_2=0$. 

\begin{proposition}\label{upperbound}
Let $v\in\cM_k$. For every $\epsilon>0$ there exist $\delta>0$ and an integer $N>0$ such that
\begin{align*}
\frac{\log\|F_{C,n}(1,s)\|_v}{\deg(F_{C,n})}-\frac{d-1}{d_1+(d-1)\deg(c)}\cdot\gamma_v<\epsilon,
\end{align*}
for all $|s|_v<\delta$ and $n\ge N$.
\end{proposition}

\begin{proof}
Let $1>\epsilon>0$. By Lemma \ref{leadingtermslimit} there exists large $N\in\N$ such that
\begin{align}\label{constanttermsmall}
\max\{|A_{C,N}(1,0)|_v,|B_{C,N}(1,0)|_v\}<(1+\epsilon/8)^{d^N}e^{\gamma_v\cdot d^N}.
\end{align}
Moreover, invoking Lemma \ref{degree}, we may choose $N\in\N$ large enough such that 
$$\frac{d^{N+i}}{\deg(F_{C,N+i})}<\frac{d-1}{d_1+(d-1)\deg(c)}+\frac{\epsilon}{8\max\{1,\gamma_v\}},$$
 for all $i\ge 0$ and we may further assume that $\frac{\log 8}{d^N}<\epsilon/16$.

Let $L=8(1+\epsilon/4)^{d^N}e^{\gamma_v\cdot d^N}$. By \eqref{constanttermsmall} we can find some $0<\delta<1$ such that for $|s|_v<\delta$ we have
\begin{align}\label{firststepfcn}
\|F_{C,N}(1,s)\|_v<\frac{L}{8}.
\end{align}
Recall from Lemma \ref{gcd} that $F_{C,n}(1,s)=F^n_{1,s}(C(1,s))$ for all $n\in\N$. From the expression of $F_{1,s}=(z^{d/2}w^{d/2}+s(\cdots), (zw)^{d_1}(z^2+w^2)+s(\cdots))$, shrinking $\delta$ if necessary and applying $F_{1,s}$ repeatedly to \eqref{firststepfcn}, we get
\begin{align*}
\|F_{C,N+i}(1,s)\|_v<\frac{L^{d^i}}{8},
\end{align*}
for all $i\ge 0$. Therefore, recalling the definition of $L$, we get
\begin{align*}
\frac{\log\|F_{C,N+i}(1,s)\|_v}{\deg(F_{C,N+i})}&<\frac{d^i}{\deg(F_{C,N+i})}\log L-\frac{\log8}{\deg(F_{C,N+i})}\\
&\le \frac{d^i}{\deg(F_{C,N+i})}\log 8+\frac{d^{N+i}}{\deg(F_{C,N+i})}\log (1+\epsilon/4)+\frac{d^{N+i}}{\deg(F_{C,N+i})}\gamma_v\\
&=\frac{d^{i+N}}{\deg(F_{C,N+i})}\cdot\frac{\log 8}{d^N}+\frac{d^{N+i}}{\deg(F_{C,N+i})}\epsilon/4+\frac{d^{N+i}}{\deg(F_{C,N+i})}\gamma_v.
\end{align*}
This inequality combined with our assumptions on $N\in\N$ yield
\begin{align*}
\frac{\log\|F_{C,N+i}(1,s)\|_v}{\deg(F_{C,N+i})}<\frac{d-1}{d_1+(d-1)\deg(c)}\cdot\gamma_v+\epsilon,
\end{align*}
 for all $i\ge 0$ and $|s|_v<\delta$. The proposition follows.
\qed

To show that the convergence is uniform from below, we will first need the following weaker estimate.
\begin{lemma}\label{weakbound}
Let $v\in\cM_k$. For all $s\in\C_v$ with $|s|_v\le 1$ and $(z,w)\in\mathbb{C}_v^2\setminus\{(0,0)\}$, we have
\begin{align*}
\frac{\|F_{1,s}(z,w)\|_v}{\|(z,w)\|^d_v}\geq |\tau|_v^{-1}\frac{|s|_v^{3\cdot d_2-2}}{4^{d-1}}.
\end{align*}
 Consequently, for all for all $n> j$ and $s\in\C_v$ with $|s|_v\leq |\tau|_v$ we have
\begin{align*}
\frac{\log\|F_{C,n}(1,s)\|_v}{d^n}\ge\frac{\log\|F_{C,j}(1,s)\|_v}{d^j}+\frac{\log(|s|_v^{3d_2-2})}{d^j}-\frac{\log(|\tau|_v4^{d-1})}{d^{j}}.
\end{align*}
\end{lemma}
\begin{proof}
Set $G_s(z,w)=(s\lambda zw,sz^2+zw+sw^2)$. We will first see that for all $s\in\C_v$ with $|s|_v\le 1$ we have
\begin{align}\label{gsinduction}
||G_s(z,w)||_v\ge\frac{|s|^2_v}{4}\|(z,w)\|^2_v.
\end{align}
Note that by the homogeniety  and symmetry of $G_s$, we may assume that $w=1$ and $|z|_v\le 1$. Then we have that either $|\lambda sz|_v\ge \frac{|s|^2_v}{4}$, or $|z|_v<\frac{|s|_v}{4}$, in which case $|s+z+sz^2|_v\ge|s|_v-\frac{|s|_v}{4}-\frac{|s|_v^3}{16}\ge\frac{1}{4}|s|^2_v$. In both cases $||G_s(z,1)||_v\ge\frac{|s|^2_v}{4}$, and our claim follows.
Note that from Proposition \ref{degenerations} and the definition of $F_{t_1, t_2}$ in \eqref{tau}, we have
$$\|F_{1,s}(z,w)\|_v=\|\tau^{-1}s^{-d_2}\cdot G^\ell_s(z,w)\|_v=|\tau|_v^{-1}|s|_v^{-d_2}\cdot\| G_s\left(G^{\ell-1}_s(z,w)\right)\|_v.$$
Repeated applications of \eqref{gsinduction} now give
\begin{align*}
\|F_{1,s}(z,w)\|_v&\ge |\tau|_v^{-1} \frac{|s|_v^{2+2^2+\cdots+2^\ell-d_2}}{4^{1+2+2^2+\cdots +2^{\ell-1}}} \|(z,w)\|_v^d=|\tau|_v^{-1} \frac{|s|_v^{3\cdot d_2-2}}{4^{d-1}} \|(z,w)\|_v^d.
\end{align*}
The first conclusion of lemma follows. For the second conclusion of the lemma, we just need to apply the first conclusion to the following and take the logarithm for both sides
  $$\left\|\frac{F_{C, n}(1,s)}{F_{C, j}(1,s)^{d^{n-j}}}\right\|_v=\left\|\frac{F_{C, n}(1,s)}{F_{C, n-1}(1,s)^{d}}\cdot \frac{F_{C, n-1}(1,s)^d}{F_{C, n-2}(1,s)^{d^{2}}}\cdots \frac{F_{C, j+1}(1,s)^{d^{n-j-1}}}{F_{C, j}(1,s)^{d^{n-j}}}\right\|_v.$$
The lemma follows. 
\end{proof}

\begin{proposition}\label{lowerbound}
Let $v\in\cM_k$. For every $\epsilon>0$, there exists $\delta>0$ and an integer $N>0$ such that
\begin{align*}
\frac{\log\|F_{C,n}(1,s)\|_v}{\deg(F_{C,n})}-\frac{d-1}{d_1+(d-1)\deg(c)}\cdot\gamma_v>-\epsilon
\end{align*}
for all $|s|_v<\delta$ and $n\ge N$.
\end{proposition}
\proof 
Recall our notation $F_{C,n}(1,s)=(A_n(s),B_n(s))=(a_n+sp_n(s),b_n+sq_n(s)).$
Let $\epsilon>0$ be small. Since by Lemma \ref{gcd} we have $F_{C, n+1}(1,s)=F_{1,s}(F_{C, n}(1,s))$, one can find
  $$\Phi_3(x,y,z,w), \Psi_3(x,y,z,w)\in k[s][x,y,z,w]$$
depending only on $\ell$ and homogenous in $\mathbf{x}:=(x,y,z,w)$ of degree $d$, such that 
\begin{equation}\label{recursivepq}
p_{n+1}(s)=\Phi_3(a_n, b_n, p_n(s), q_n(s))\textup{ and } q_{n+1}=\Psi_3(a_n, b_n, p_n(s), q_n(s)),
\end{equation}
for all $n\geq 0$. Moreover, one can find a large $L_0>0$ such that 
 \begin{equation}\label{qpgrowth}\|\Phi_3(\mathbf{x})\|_v\leq L_0 \cdot \|\mathbf{x}\|_v^d \textup{ and } \|\Psi_3(\mathbf{x})\|_v\leq L_0 \cdot \|\mathbf{x}\|_v^d.
 \end{equation}
Enlarging $L_0$ if necessary, we may assume 
\begin{equation}\label{boundcepsilons}
(3d_2-2)\cdot d\cdot\log\left(\frac{1-\epsilon/L_0}{1+\epsilon/L_0}\right)> -\frac{\epsilon}{4}.
\end{equation}In view of Lemma \ref{leadingtermslimit} we can find a large $N\in\N$  such that
\begin{align}\label{preconstant}
\max\{|p_{N}(0)|_v,|q_{N}(0)|_v\}<\left(1+\frac{\epsilon}{L_0}\right)^{d^{N}}\cdot e^{\gamma_v \cdot d^{N}},
\end{align}
and also
\begin{align}\label{constantbig}
\left(1-\frac{\epsilon}{L_0}\right)^{d^{N+i}}\cdot e^{\gamma_v \cdot d^{N+i}}<|a_{N+i}|_v,|b_{N+i}|_v<\left(1+\frac{\epsilon}{L_0}\right)^{d^{N+i}}\cdot e^{\gamma_v \cdot d^{N+i}},
\end{align}
for all $i\ge 0$. By Lemma \ref{degree}, enlarging $N$ if necessary, we may further assume that 
\begin{align}\label{tofixthedegree}
\frac{d^{N+i}}{\deg(F_{C,N+i})}>\frac{d-1}{d_1+(d-1)\deg(c)}-\frac{\epsilon}{L_0\cdot \max\{1,\gamma_v\}},
\end{align}
for all $i\ge 0$ and that
\begin{align}\label{boundcepsilon}
(3d_2-2)\left(\frac{\log\epsilon}{d^{j}}-\frac{\log L_0}{d^{N-1}}\right) -\frac{\log(|\tau|_v4^{d-1})}{d^{j}}> -\frac{\epsilon}{4}
\end{align}
for all $j\geq N$. 

Define $L:=L_0\cdot \left((1+\epsilon/L_0)e^{\gamma_v}\right)^{d^N}$.
By \eqref{preconstant} and \eqref{constantbig}, we can find a small $0<\delta<\min\{1,|\tau|_v\}$ such that if $|s|_v<\delta$, then
\begin{align}
\max\{|p_N(s)|_v,|q_N(s)|_v\}<\frac{L}{L_0}.
\end{align}
Combining this with the recursive relations defining $p_n(s), q_n(s)$ given in \eqref{recursivepq} and inequalities \eqref{qpgrowth} and \eqref{constantbig}, we get inductively that if $|s|<\delta$, then for all $i\ge 0$ we have
\begin{align}\label{uglyformulainuse}
\max\{|p_{N+i}(s)|_v, |q_{N+i}(s)|_v\}<\frac{L^{d^i}}{L_0}.
\end{align}
Now choose an integer $N'>N$ such that
\begin{align*}
\delta':=\frac{\epsilon(1-\epsilon/L_0)^{d^{N'}}}{(1+\epsilon/L_0)^{d^{N'}}L_0^{d^{N'-N}}}<\delta.
\end{align*}
We will show that if $|s|_v<\delta'$ and $n\ge N'$ we have 
$$\frac{\log\|F_{C,n}(1,s)\|_v}{\deg(F_{C,n})}\ge \frac{d-1}{d_1+(d-1)\deg(c)} -\epsilon.$$
Denote by $\delta_j:=\frac{\epsilon(1-\epsilon/L_0)^{d^{j}}}{(1+\epsilon/L_0)^{d^{j}}L_0^{d^{j-N}}}>0$.        
If $j\ge N'>N$ and $|s|_v\le\delta_j\le\delta'<\delta$ we have 
\begin{align}\label{fornandn}
\begin{split}
\frac{\log\|F_{C,j}(1,s)\|_v}{d^j}&\ge \frac{\log|a_j+sp_j(s)|_v}{d^j}\ge \frac{\log(|a_j|_v-|sp_j(s)|_v)}{d^j}
\\&\ge \frac{\left(1-\frac{\epsilon}{L_0}\right)^{d^{j}}\cdot e^{\gamma_v \cdot d^{j}}-\frac{|s|_vL^{j-N}}{L_0}}{d^j}, \text{ by  \eqref{constantbig} and \eqref{uglyformulainuse}}.
\end{split}
\end{align}
Therefore,
\begin{align}\label{fornandn2}
\begin{split}
\frac{\log\|F_{C,j}(1,s)\|_v}{d^j}
&\ge\log\left(1-\frac{\epsilon}{L_0}\right)+ \gamma_v+\log\left(1-\frac{\epsilon}{L_0}\right) \ge -\frac{\epsilon}{4} +\gamma_v,
\end{split}
\end{align}
when $\epsilon>0$ is sufficiently small and $L_0$ is large enough. 
Combining this with \eqref{tofixthedegree} we get
\begin{align*}
\frac{\log\|F_{C,j}(1,s)\|_v}{\deg(F_{C,j})}\ge \frac{d-1}{d_1+(d-1)\deg(c)}\gamma_v-\epsilon.
\end{align*}
Thus if for $n\ge N'$ we have $|s|_v\le\delta_n<\delta'$, then
\begin{align*}
\frac{\log\|F_{C,n}(1,s)\|_v}{\deg(F_{C,n})}\ge\frac{d-1}{d_1+(d-1)\deg(c)}\gamma_v-\epsilon.
\end{align*}
On the other hand, if $n\ge N'$ and $\delta_n< |s|_v<\delta'$,  then there is some $j$ with $N'\le j<n$ such that $\delta_{j+1}< |s|_v\le\delta_{j}<|\tau|_v$.
By Lemma \ref{weakbound} we have
\begin{align*}
\frac{\log\|F_{C,n}(1,s)\|_v}{d^n}\ge\frac{\log\|F_{C,j}(1,s)\|_v}{d^j}+\frac{\log(|s|_v^{3d_2-2})}{d^j}-\frac{(d-1)\log4}{d^{j}}.
\end{align*}
This upon using \eqref{boundcepsilon}  and \eqref{fornandn2} yields
\begin{align*}
\frac{\log\|F_{C,n}(1,s)\|_v}{d^n}&\ge -\frac{\epsilon}{4} +\gamma_v+\frac{\log(|s|_v^{3d_2-2})}{d^j}-\frac{(d-1)\log4}{d^{j}}\\
&\ge-\frac{\epsilon}{4} +\gamma_v+(3d_2-2)\left(\frac{\log\epsilon}{d^{j}}-\frac{\log L_0}{d^{N-1}}\right) +d\cdot\log\left(\frac{1-\epsilon/L_0}{1+\epsilon/L_0}\right)-\frac{\log(|\tau|_v4^{d-1})}{d^{j}}\\
&\ge -\frac{\epsilon}{4} +\gamma_v+(3d_2-2)\cdot d\cdot\log\left(\frac{1-\epsilon/L_0}{1+\epsilon/L_0}\right)-\frac{\epsilon}{4} \text{, by \eqref{boundcepsilons}}\\
&\ge\gamma_v-\frac{\epsilon}{4} -\frac{\epsilon}{4}-\frac{\epsilon}{4}= \gamma_v-\frac{3\epsilon}{4},
\end{align*}
for sufficiently small $\epsilon>0$. Finally, upon using \eqref{tofixthedegree} we get
\begin{align*}
\frac{\log\|F_{C,n}(1,s)\|_v}{\deg(F_{C,n})}\ge\frac{d-1}{d_1+(d-1)\deg(c)}\gamma_v-\epsilon,
\end{align*}
as claimed. This finishes the proof.
\end{proof}

\emph{Proof of Theorem \ref{continuity}.}
By a standard telescoping sum argument as in  \cite[Proposition 1.2]{Branner:Hubbard:1} we see that the functions $\frac{\log\|F_{C,n}(t_1,t_2)\|_v}{\deg(F_{C,n})}$ converge locally uniformly to the function $G_{F,C,v}$ on $\C_v^2\setminus \C_v\times \{0\}$.
Thus, it suffices to prove that the sequence 
   $$\frac{\log\|F_{C,n}(1,s)\|_v}{\deg(F_{C,n})}$$
converges locally uniformly in a neighborhood of $s=0$. This now follows from Propositions \ref{upperbound} and \ref{lowerbound}. \qed
\subsection{Bounds of the radii}
Recall that we have chosen the lift of $c\in k(t)$, denoted by $C(t_1,t_2)=(A(t_1,t_2),B(t_1,t_2))$, so that the coefficients of $C(t_1,t_2)$ lie in $\mathcal{O}_k$. In particular, $\Res(A,B)\in\mathcal{O}_k$. In what follows, we let
\begin{itemize}
\item
$S\subset\cM_{k}$ be the finite set consisting of the non-archimedean places of $k$ such that if $v\in S$:
$$|\Res(A,B)|_v<1\textup{ or } |\tau|_v<1 \textup{ or } |\Res_{(z,w)}(F_{t,1})|_v<1,$$
where $\tau\in\Z[\la]\subset\mathcal{O}_k$ is the one defined in Proposition \ref{degenerations}.
\item
$\cM'_{k}$ be the set of non-archimedean places in $\cM_k$ satisfying the following. If  $v\in\cM'_k$, then $|\tau|_v=1$ and there exists $n_v\in\N$ such that $$\|F_{C,n_v}(1,0)\|_v=\max\{|A_{C,n_v}(1,0)|_v,|B_{C,n_v}(1,0)|_v\}<1.$$
\item
$\cM'_{k,m}$ the set of places $v\in\cM'_k$ for which
$m\in\N$ is the smallest integer such that $$|B_{C,m}(1,0)|_v<1.$$
\end{itemize}
Furthermore, we denote the non-archimedean places of $k$ by $\cM^0_k$, and the archimedean ones by $\cM^{\infty}_k$.
Finally, note that from equation \eqref{anbnourpair} we get
 $$\cM'_k=\displaystyle\bigcup_{\substack{m\in\N}}\cM'_{k,m}.$$
\begin{lemma}\label{smallteasy}
Let $v\in\cM^0_{k}\setminus S$. If $t\in\C_v$ has $|t|_v\le 1$, then $\|F_{C,n}(t,1)\|_v=1$ for all $n\in\N$. In particular,
$$G_{F, C,v}(t,1)=0.$$
\end{lemma}
\begin{proof}
Let $v\in\cM^0_{k}\setminus S$ and consider $t\in\C_v$ with $|t|_v\le 1$. Since the coefficients of $C(t_1,t_2)$ are $v-$adic integers, we conclude $|\Res(A,B)|_v=\|(t,1)\|_v=1$. 
Then \cite[Lemma 10.1]{BRbook} yields  $\|C(t,1)\|_v=1$. 
Since $v\notin S$ and $\tau\in\mathcal{O}_k$, we have $|\tau|_v=1$. Now from the definition of $F_{t_1,t_2}$ in \eqref{tau}, we see that all its coefficients are $v$-adic integers. 
As furthermore $|\Res_{(z,w)}(F_{t,1})|_v=1$, using \cite[Lemma 10.1]{BRbook} once more, we get $\|F^n_{t,1}(C(t,1))\|_v=1$ for all $n\in\N$. This in turn, by Lemma \ref{gcd}, implies $\|F_{C,n}(t,1)\|_v=1$ for all $n\in\N$ as claimed. 
\end{proof}
We write $M_{C,v}:=\{(t_1,t_2)\in\C_v^2\setminus\{(0,0)\}~:~G_{F,C,v}(t_1,t_2)\le 0\}$. 
\begin{proposition}\label{trivialsets} 
For all places $v\in\cM^0_k\setminus(S\cup\cM'_k)$, we have
$$G_{F, C,v}(t_1,t_2)=\log\|(t_1,t_2)\|_v.$$
In particular, $r_{in}(M_{C,v})=r_{out}(M_{C,v})=1$.
\end{proposition}
\begin{proof}
Let $v\in\cM^0_k\setminus(\cM'_k\cup S)$ and $(t_1,t_2)\in\C_v^2\setminus\{(0,0)\}$. Since $G_{F, C,v}$ scales logarithmically, by Lemma \ref{smallteasy}, we know that the claim holds when $|t_1|_v\le|t_2|_v$. It suffices to show that for $t_2\in\C_v$ with $|t_2|_v<1$, we have $\|F_{C,n}(1,t_2)\|_v=1$ for all $n\in\N$. 

Assume to the contrary that there exist some $t_2\in\C_v$ with $|t_2|_v<1$ and  $n\in\N$ such that $\|F_{C,n}(1,t_2)\|_v\neq1.$
Notice that since $v\notin S$, we have $|\tau|_v=1$, which using equation \eqref{tau} and Lemma \ref{gcd}, yields that $F_{C,n}$ has integral coefficients. Hence our assumption implies 
$$\|F_{C, n}(1,t_2)\|_v=\|(A_{C, n}(1, 0)+t_2(\cdots), B_{C, n}(1, 0)+t_2(\cdots))\|_v< 1.$$
As we also have $|t_2|_v<1$, this gives
$$\|(A_{C, n}(1, 0), B_{C, n}(1, 0))\|_v<1.$$
Since also $|\tau|_v=1$, we have $v\in\cM'_k$, which is a contradiction. This finishes our proof.
\end{proof}
\begin{lemma}\label{forlowerbound}
Let $m\in\N$ and $v\in\cM'_{k,m}\setminus S$. For all $n\ge m+1$ we have
\begin{align*}
|A_{C,n}(1,0)|_v<|B_{C,n}(1,0)|_v.
\end{align*}
In particular, for all $n\ge m+1$ we have
$$|A_{C,n}(1,0)|_v^d<| A_{C,n+1}(1,0)|_v, ~|B_{C,n+1}(1,0)|_v<|B_{C,n}(1,0)|_v^d.$$
Moreover, $\displaystyle\left\{d^{-n}\log|A_{C,n}(1,0)|_v\right\}_{n\ge m+1}$ is increasing and
$\displaystyle\left\{d^{-n}\log|B_{C,n}(1,0)|_v\right\}_{n\ge m+1}$ is decreasing.
\end{lemma}
\begin{proof}
Let $m\in\N$ and $v\in \cM'_{k,m}\setminus S$. Recall that from \eqref{anbnourpair}, we have
\begin{align*}
a_{n+1}=a_n^{d_1}b_n^{d_1}\cdot a_nb_n, ~b_{n+1}=a^{d_1}_nb^{d_1}_n(a_n^2+b_n^2).
\end{align*}
Since $v\in\cM'_{k,m}$ we have $|b_m|_v<1$ and $|b_n|_v=1$ for all $n<m$. We will prove the first inequality in the lemma by induction. Using the ultrametric inequality, it is easy to see that if for some $n\in\N$ we have $|a_n|_v<|b_n|_v$ then $|a_{n+1}|_v<|b_{n+1}|_v$. Thus it remains to prove the base case; that is $|a_{m+1}|_v<|b_{m+1}|_v$. To this end, we consider cases depending on the value of $m\in\N$. 

If $m=1$, we have $|b_1|_v<1$ and $|a_0|_v\le |b_0|_v=1$. To see that $|a_2|_v=|a_1b_1|_v^{d_1+1}<|b_2|_v=|(a_1b_1)^{d_1}(a_1^2+b_1^2)|_v$, it suffices to show $|a_1b_1|_v<|a_1^2+b_1^2|_v$. Note that $|a_1|_v=|a_0|^{d_1+1}_v$ and $|b_1|_v=|a_0|_v^{d_1}|a_0^2+b_0^2|_v$. If $|a_0|_v<1$ we have $|a_1|_v<|b_1|_v$; hence $|a_1b_1|_v<|a_1^2+b_1^2|_v$. If $|a_0|_v=1$, then $|a_1|_v=1$ and  $|a_1b_1|_v=|b_1|_v<|a_1^2+b_1^2|_v=1$ holds as well, since $|b_1|_v<1$.
The base case follows.

If on the other hand $m=0$ or $m\ge 2$, we will prove that $|a_{m}|_v=1$. Then $|a_{m+1}|_v=|b_m|_v^{d_2}<|b_m|_v^{d_1}=|b_{m+1}|_v$ and the base case follows. If $m=0$, so that $|b_0|_v<1$, our assumption that $v\notin S$ and hence $|\Res(A,B)|_v=1$ yields $|a_0|_v=1$, as claimed. If now $m\ge 2$, assume that $|a_{m}|_v<1$ to end in a contradiction. Since $v\in\cM'_{k,m}$ we have $|b_{m-1}|_v=|b_{m-2}|=1$, which by \eqref{anbnourpair} and our assumption that $|a_m|_v<1$ implies $|a_{m-1}|_v<1$ and $|a_{m-2}|_v<1$. This in turn gives $|b_{m-1}|_v<1$ contradicting the minimality of $m\in\N$. Hence $|a_{m}|_v=1$ and the base case follows as noted.
This completes the proof of the first inequality in the lemma. The rest of the lemma now follows by applying \eqref{anbnourpair}.
\end{proof}

In the course of the proof of Lemma \ref{forlowerbound}, we saw the following equalities, which will be handy later on. They are an immediate consequence of the ultrametric inequality, the definition of $\cM'_{k,m}$ and the recursive definition of $a_n$ and $b_n$ in \eqref{anbnourpair}. 
\begin{remark}
Let $m\in\N_{\ge 2}$ and $v\in\cM'_{k,m}$. We have
\begin{align}\label{amamminusone}
|A_{C,m}(1,0)|_v=|A_{C,m-1}(1,0)|_v=|B_{C,m-1}(1,0)|_v=1.
\end{align}
In particular,
\begin{align}\label{fromatobwhenmbig}
|A_{C,m+1}(1,0)|_v=|B_{C,m}(1,0)|^{d_2}_v=|B_{C,m+1}(1,0)|^{d_2/d_1}_v.
\end{align}
\end{remark}
Before stating the next proposition, recall that from Lemma \ref{degree} we have  $\hat{h}_f(c)\neq 0$.
\begin{proposition}\label{radiibound}
Let $m\in\N_{\ge 2}$ and $v\in\cM'_{k,m}\setminus S$. We have $$\bar{D}^2(0,1)\subset M_{C,v}\subset \bar{D}^2\left(0,e^{-\frac{1}{\hat{h}_f(c)}\left( \frac{3d_2\cdot \log|a_{m+1}|_v}{d^{m}}-\frac{\log4^{d-1}}{d^{m}}\right)}\right).$$
In particular,
\begin{align*}
 0\le \log r_{in}(M_{C,v}) \le \log r_{out}(M_{C,v})\le -\frac{1}{\hat{h}_f(c)}\cdot\left( \frac{3d_2\cdot \log|a_{m+1}|_v}{d^{m}}-\frac{\log4^{d-1}}{d^{m}}\right).
\end{align*}
\end{proposition}

\begin{proof}
Let $m\ge 2$ and $v\in\cM'_{k,m}\setminus S$. Since the coefficients of $F_{C,n}$ are $v$-adic integers for all $n\in\N$, we have $\|F_{C,n}(t_1,t_2)\|_v\le\|(t_1,t_2)\|_v^{\deg(F_{C,n})}$. Therefore $\bar{D}^2(0,1)\subset M_{C,v}$ and the first inclusion follows.
By Lemma \ref{smallteasy} we get that if $|t|_v\le 1$, then $G_{F,C,v}(t,1)=0$. Therefore to show the reverse inclusion it suffices to prove that if $0<|s|_v<1$, then 
\begin{align}\label{innerradboundfc}
G_{F, C,v}(1,s)\ge\frac{1}{\hat{h}_f(c)}\cdot\left( \frac{3d_2\cdot \log|a_{m+1}|_v}{d^{m}}-\frac{\log4^{d-1}}{d^{m}}\right).
\end{align}
To this end, let $s\in\mathbb{C}_v$ be such that $0<|s|_v<1$. Recall from Lemma \ref{forlowerbound} that $\{|b_n|_v\}_{n\ge m+1}$ is a strictly decreasing sequence which converges to zero. 
Assume first that $0<|s|_v<|b_{m+1}|_v$. Then there is some $j\ge m+1$ such that $|b_{j+1}|_v\le |s|_v<|b_j|_v$.  As the coefficients of $F_{C,j}(1,s)=(a_{j}+s(\cdots),b_{j}+s(\cdots))$ are $v$-adic integers and since from Lemma \ref{forlowerbound} we have $|a_j|_v<|b_j|_v$ for all $j\ge m+1$, we get $\|F_{C,j}(1,s)\|_v=|b_j|_v$. Therefore, upon using Lemma \ref{weakbound} (note that here $|\tau|_v=1$) we get
\begin{align*}
\frac{\log\|F_{C,n}(1,s)\|_v}{d^n}&\ge\frac{\log\|F_{C,j}(1,s)\|_v}{d^j}+\frac{\log(|s|_v^{3d_2-2})}{d^j}-\frac{\log(|\tau|_v4^{d-1})}{d^{j}}\\
&\ge\frac{\log|b_j|_v}{d^j}+d(3d_2-2)\frac{\log|b_{j+1}|_v}{d^{j+1}}-\frac{\log4^{d-1}}{d^{j}},
\end{align*}
for all $n>j$. This in turn, using Lemma  \ref{forlowerbound} implies
\begin{align}\label{case2mstep2}
\begin{split}
\frac{\log\|F_{C,n}(1,s)\|_v}{d^n}&\ge(d(3d_2-2)+1)\frac{\log|b_{j+1}|_v}{d^{j+1}}-\frac{\log4^{d-1}}{d^{j}}\\
&\ge (d(3d_2-2)+1)\frac{\log|a_{j+1}|_v}{d^{j+1}}-\frac{\log4^{d-1}}{d^{j}}\\
&\ge(d(3d_2-2)+1)\frac{\log|a_{m+1}|_v}{d^{m+1}}-\frac{\log4^{d-1}}{d^{m+1}}\\
&\geq  3d_2\frac{\log|a_{m+1}|_v}{d^{m}}-\frac{\log4^{d-1}}{d^{m}}.
\end{split}
\end{align}
Finally, assume that $|b_{m+1}|_v\leq |s|_v<1$. 
From equation \eqref{amamminusone} we have $|b_m|_v<|a_m|_v=1$; thus $\|F_{C, m}(1,s)\|_v=1$. Invoking now Lemmata \ref{weakbound} and \ref{forlowerbound}, we get
\begin{align}\label{case1m}
\begin{split}
\frac{\log\|F_{C,n}(1,s)\|_v}{d^n}&\ge
\frac{\log(|s|_v^{3d_2-2})}{d^m}-\frac{\log4^{d-1}}{d^{m}}\ge\frac{\log(|b_{m+1}|_v^{3d_2-2})}{d^m}-\frac{\log4^{d-1}}{d^{m}}\\
&\geq \frac{\log(|a_{m+1}|_v^{3d_2-2})}{d^m}-\frac{\log4^{d-1}}{d^{m}}\geq  3d_2\frac{ \log|a_{m+1}|_v}{d^{m}}-\frac{\log4^{d-1}}{d^{m}}, 
\end{split}
\end{align}
for all $n>m$. Letting $n\to\infty$ in \eqref{case2mstep2} and  \eqref{case1m}, inequality \eqref{innerradboundfc} follows. This finishes our proof.
\end{proof}

\begin{lemma}\label{numberofprimes}
There exist constants $L_1,L_2>0$ such that ]$\sum_{v\in \cM_{k, m}\backslash S} N_v\le L_1\cdot2^m$  and $\displaystyle\prod_{\substack{v\in\cM'_{k,m}}}|A_{C,m+1}(1,0)|^{-N_v}_v<L_2^{2^m}$ for all $m\ge 2$. 
\end{lemma}
\begin{proof}
We write  $T_m:=\sum_{v\in \cM_{k, m}\backslash S} N_v$.
Firstly we are going to prove that for some $L_1>0$ we have $T_m<L_1\cdot 2^m$. To this end, we define
$\mathcal{P}_m:=\{v\in\cM^0_k\setminus S~:~ |b^*_{m}|_v<1\}$. We claim that $\cM'_{k, m}\backslash S\subset \mathcal P_m$; hence it suffices to prove that for some $L_1>0$ we have $\sum_{v\in \mathcal P_m} N_v<L_1\cdot 2^m$. To prove that our claim holds, let $m\ge2$ and $v\in \cM'_{k, m}\backslash S$. We recall from \eqref{amamminusone} that $|a_m|_v=1$. Moreover, from \eqref{anbnstar} and \eqref{aldef}, we have that $a_m$ and $a^*_m$ are $v-$adic integers and $|a_m|_v=|a_m^*|_v|\al_m|_v=1$. Hence $|a_m^*|_v=|\al_m|_v=1$ and 
\begin{align}\label{bmbmstar}
|b_m|_v=|\al_m|_v\cdot |b^*_m|_v=|b^*_m|_v.
\end{align}
Our claim follows. Now notice that the recursive definition of $\{b^*_n\}$ in \eqref{anbnstar} allows us to conclude that there is a constant $L_0>1$ such that 
\begin{align}
\label{easyuper}
\prod_{\substack{v\in\cM_k~ : ~|b^*_{m}|_v>1}}|b^*_{m}|_v^{N_v}\le L_0^{2^m}.
\end{align}
Let $r:=\displaystyle\sup_{v\in \cM^0_k}\{ |\alpha|_v: |\alpha|_v<1 \textup{ and } \alpha \in k\}\in(0,1)$. For each $v\in\mathcal{P}_m$ we have
\begin{align}
\label{pbiggerthantwo}
|b^*_{m}|_v^{N_v}\le r^{N_v}.
\end{align}
Combining \eqref{easyuper} and \eqref{pbiggerthantwo} and upon using the product formula we have
\begin{equation}\label{toboundplaces}
\prod_{\substack{v\in\mathcal{P}_m}} r^{N_v}\ge\prod_{\substack{v\in\mathcal{P}_m}}|b^*_{m}|_v^{N_v}\ge \prod_{\substack{v\in\cM_k ~:~ |b^*_{m}|_v<1}}|b^*_{m}|_v^{N_v}  =\prod_{\substack{v\in\cM_k ~:~ |b^*_{m}|_v>1}}|b^*_{m}|_v^{-N_v} \ge L_0^{-2^m}.
\end{equation}
Thus if $L_1=\log_{1/r} L_0$, we get $T_m\le L_1\cdot2^m$ and the first part of the lemma follows.
For the second part of the lemma, first note that \eqref{fromatobwhenmbig} combined with \eqref{bmbmstar} implies
\begin{align}\label{nostarstart}
|a_{m+1}|_v=|b_m|^{d_2}_v=|b^*_m|_v^{d_2}.
\end{align}
Using this and the fact that $\cM'_{k, m}\backslash S\subset \mathcal P_m$, inequality \eqref{toboundplaces} yields
$$
\displaystyle\prod_{\substack{v\in\cM'_{k,m}\backslash S}}|a_{m+1}|^{-N_v}_v=\displaystyle\prod_{\substack{v\in\cM'_{k,m}\backslash S}}|b^*_{m}|^{-{d_2}\cdot N_v}_v
\le\displaystyle\prod_{\substack{v\in\mathcal{P}_m}}|b^*_{m}|^{-{d_2}\cdot N_v}_v \le L_0^{d_2\cdot 2^{m}}.$$
Setting $L_2=L_0^{d_2}$, the lemma follows.
\end{proof}
We can now control the products of the inner and outer radii as in the following proposition.
\begin{proposition}\label{strongconvergence}
The products
\begin{align*}
\prod_{v\in \cM_k} r_{out}(M_{C,v})^{N_v}\text{, }\prod_{v\in \cM_k} r_{in}(M_{C,v})^{N_v}\text{ and }  \prod_{v\in \cM_k}\capacity(M_{C,v})^{N_v}
\end{align*}
converge strongly.
\end{proposition}
\begin{proof}
We first note that the set $\cM_k^{\infty}\cup S\cup\cM'_{k,0}\cup\cM'_{k,1}$ is finite. The continuity of the potentials, proved in Theorem \ref{continuity}, yields that the products 
\begin{align*}
\prod_{v\in\cM_k^{\infty}\cup S\cup\cM'_{k,0}\cup\cM'_{k,1}} r_{out}(M_{C,v})^{N_v}\text{ and }\prod_{v\in\cM_k^{\infty}\cup S\cup\cM'_{k,0}\cup\cM'_{k,1}} r_{in}(M_{C,v})^{N_v}
\end{align*}
are finite. Hence, invoking Theorem \ref{quasi-adelic set}, Proposition \ref{trivialsets} and Proposition \ref{radiibound}, it suffices to prove that the following sum converges
\begin{align*}
\sum_{m=2}^{\infty}\sum_{\substack{v\in\cM'_{k,m}}}N_v\cdot \left(-\frac{3d_2\cdot \log|a_{m+1}|_v}{d^{m}}+\frac{\log4^{d-1}}{d^{m}}\right).
\end{align*}
This in turn follows from Lemma \ref{numberofprimes}. 
\end{proof}
\subsection{Proof of Theorem \ref{mainthm3}}
Let $k$ be a number field and $c\in k(t)$ be such that $0,\infty\notin\mathcal{O}_{f_{\infty}}(c(\infty))$. As $f=g^{\ell}_{\la}$, it suffices to prove the conclusions of our theorem for the pair $(f,c)$ in place of $(g_{\la},c)$; see Remark \ref{moregeneric}.
First we are going to see that the measure $\mu_{f,c}$ is quasi-adelic. Since for each $v\in\cM_k$ we have $\mu_{f,c,v}=\mu_{M_{C,v}}$, by Theorem \ref{quasi-adelic set} it suffices to prove that the set $\{M_{C,v}\}_{v\in\cM_k}$ is quasi-adelic. The continuity of the potentials of $M_{C,v}$ is established in Theorem \ref{continuity}. Moreover, in view of Proposition \ref{strongconvergence} we know that the products
$$\prod_{v\in \cM_k} r_{out}(M_{C,v})^{N_v}\text{ and }\prod_{v\in \cM_k} r_{in}(M_{C,v})^{N_v}$$
converge strongly. 
Thus the measure $\mu_{f,c}$ is quasi-adelic. Assume further that $c(\infty)$ is not a preperiodic point for $f_{\infty}$. Then our assumption on $c$ implies that for $\infty\in\Sing(f)$, the dynamical pair $(g,c)$ is $\infty-$generic. By Theorem \ref{genericnoadelic} we conclude that $(f,c)$ is not adelic. 
\qed

\section{Variation of canonical heights and equidistribution on $\P^1$}\label{Svchunlikely}
Our aim in this section is to prove Theorem \ref{vchunlikely}.
In what follows, we keep the notation as in Section \ref{a pair is rarely adelic}. In what follows we let $F$ and $C=(A,B)$ of $f$ and $c$ respectively and enlarge the number field $k$ if necessary so that $F$ and $C$ are defined over $k$ and $\Sing(f)\subset k$. We denote by
$$G_{F,C,v}(t_1,t_2)=\lim_{n\to\infty}\frac{\log||F_{C,n}(t_1,t_2)||_v}{\deg F_{C,n}}.$$
We also write $M_{F,C,v}=\{(t_1,t_2)\in\C_v^2~:~G_{F,C,v}(t_1,t_2)\le 0\}.$

In the following proposition we show that the height associated with a measure $\mu_{f,c}$ for a quasi-adelic pair $(f,c)$, is proportional with the Call-Silverman canonical height; hence both heights have the same small points. The first author of this article is indebted to Laura DeMarco for many ideas in this proof.
\begin{proposition}\label{samesmallpoints}
Let $k$ be a number field and let $f\in k(t)(z)$ and $c\in k(t)$ be such that the dynamical pair $(f,c)$ is quasi-adelic. For any  $t\in\overline{k}\setminus\Sing(f)$ we have 
\begin{align*}
\hhat_{\mu_{f,c}}(t)=\frac{[k:\Q]}{\hat{h}_f(c)}\cdot\hhat_{f_{t}}(c(t)).
\end{align*}
\end{proposition}
\begin{proof}
We write $d=\deg_z f$. Let $t\in\overline{k}\setminus\Sing(f)$ and write $S=\Gal(\overline{k}/k)\cdot t$.  The definition of our height in \eqref{heightq}  gives
\begin{align*}
 \hhat_{\mu_{f,c}}(t)= \frac{1}{|S|}  \cdot \sum_{x\in S}\sum_{v\in \cM_k} \left(N_v\cdot  G_{F,C,v}(x,1)+\frac{1}{2}\log \capacity(M_{F,C,v})^{N_v}\right).
\end{align*}
First we will see that
\begin{align}\label{height0}
\frac{1}{|S|}  \cdot \sum_{x\in S}\sum_{v\in \cM_k} N_v\cdot  G_{F,C,v}(x,1)=\frac{[k:\Q]}{\hat{h}_f(c)}\cdot \hhat_{f_{t}}(c(t)).
\end{align}
To this end, notice that from the definition of the Call-Silverman canonical height, we have
\begin{align}\label{height1}
\hhat_{f_{t}}(c(t))&=\frac{1}{[k(t):\Q]}\displaystyle\lim_{n\to\infty}\sum_{\substack{x\in S}}\sum_{v\in\cM_k}N_v\frac{\log||F_{C,n}(x,1)||_v}{d^n}.
\end{align}
Arguing as in Lemma \ref{smallteasy}, we see that for all but finitely many places $v\in\cM_k$, 
we have $||F_{C,n}(x,1)||_v=1$ for all $n\in\N$. More specifically, for fixed $x\in S$, this conclusion holds for all places $v\in\cM^0_k$ such that the coefficients of $F$ and $C$ are $v-$adic integers and $|x|_v=|\Res(A,B)|_v=|\Res_{(z,w)}(F_{t,1})|_v=|u_{\beta}(x,1)|_v=1$ for all $\beta\in\Sing(f)$.  
Therefore we can interchange the limit with the summation in \eqref{height1} to get
\begin{align*}
\hhat_{f_{t}}(c(t))&=\frac{1}{[k(t):\Q]}\sum_{\substack{x\in S}}\sum_{v\in\cM_k}N_v\displaystyle\lim_{n\to\infty}\frac{\log||F_{C,n}(x,1)||_v}{d^n}\\
&=\frac{\hat{h}_f(c)}{[k(t):\Q]}\sum_{\substack{x\in S}}\sum_{v\in\cM_k}N_v\cdot G_{F,C,v}(x,1)\\
&=\frac{\hat{h}_f(c)}{[k:\Q]}\cdot\frac{1}{|S|}  \cdot \sum_{x\in S}\sum_{v\in \cM_k} N_v \cdot G_{F,C,v}(x,1).
\end{align*}
Thus we have established equation \eqref{height0}.
We now have 
\begin{align}\label{reducetoglobalcapacity}
 \hhat_{\mu_{f,c}}(t)=\frac{[k:\Q]}{\hat{h}_f(c)}\cdot \hhat_{f_{t}}(c(t))+\frac{1}{2}\log \prod_{v\in \cM_k}\capacity(M_{F,C,v})^{N_v}.
\end{align}
Since $\capacity(M_{F,C,v})^{N_v}\geq 1$ for all but finitely many $v\in \cM_k$, we see that  $\prod_{v\in \cM_k}\capacity(M_{F,C,v})^{N_v}$ converges strongly. 
It remains to show that the global logarithmic capacity is equal to zero, that is
\begin{align}\label{globalcapacity}  
\log \prod_{v\in \cM_k}\capacity(M_{F,C,v})^{N_v}=0.
\end{align}

The authors in \cite{DWY:Lattes, DWY:Per1} put a lot of effort into computing the explicit resultant formula for each $F_{C, n}$ to show that the global logarithmic capacity is zero. It is much harder to compute the resultants of $F_{C, n}$ here. Instead, we take a different approach, making use of Proposition \ref{finiteness}. Laura DeMarco has independently communicated  a similar idea with the first author of this article. 

Towards the proof of \eqref{globalcapacity}, we first note that there are infinitely many $t\in\overline{k}$ such that $\hhat_{f_t}(c(t))=0$; see \cite[Theorem 1.6]{DeMarco: heights}. 
Hence from \eqref{reducetoglobalcapacity} we get that for infinitely many $t\in\overline{k}$ we have
 $$\hhat_{\mu_{f,c}}(t)=\frac{1}{2}\log \prod_{v\in \cM_k}\capacity(M_{C,v})^{N_v}.$$
 This in turn combined with Proposition (\ref{finiteness}) yields
 $$\log \prod_{v\in \cM_k}\capacity(M_{F,C,v})^{N_v}\geq 0.$$
We have reduced our claim to proving
\begin{align}\label{negativeglobalcap}
   \log \prod_{v\in \cM_k}\capacity(M_{F,C,v})^{N_v}\leq 0.
\end{align}
From Proposition \ref{capacity resultant}, we have 
  $$\capacity(M_{F,C,v})\leq \liminf_{n\to \infty} |\Res F_{C, n}|_v^{-1/\deg(F_{C,n})^2}.$$
Note now that there exist a finite subset $S_0\subset \cM_k $ containing all archimedean places of $k$ such that the coefficients of $F_{t_1, t_2}$ and $C(t_1, t_2)$ are $S_0$-integers and the elements of $\Sing(f)$ are $S_0-$units. 
Then we have 
\begin{align}\label{intres}
N_v\cdot\frac{\log|\Res(F_{C, n})|_v}{\deg(F_{C,n})^2}\le0,
\end{align}
for all $n\in\N$ and $v\in \cM_k\backslash S_0$. Moreover, 
 \begin{align}\label{globalcapacitynonarch}
 \log \prod_{v\in \cM_k}\capacity(M_{F,C,v})^{N_v}&=\sum_{v\in \cM_k\backslash S_0}\log \capacity(M_{F,C,v})^{N_v}+\sum_{v\in S_0}\log \capacity(M_{F,C,v})^{N_v}\nonumber\\
&\leq \sum_{v\in \cM_k\backslash S_0}\log \capacity(M_{F,C,v})^{N_v}+\liminf_{n\to \infty} \sum_{v\in S_0}-N_v\cdot\frac{\log|\Res(F_{C, n})|_v}{\deg(F_{C,n})^2}\\
    &=\sum_{v\in \cM_k\backslash S_0}\log \capacity(M_{F,C,v})^{N_v}+\liminf_{n\to \infty}\sum_{v\in \cM_k\backslash S_0}N_v\cdot\frac{\log|\Res(F_{C, n})|_v}{\deg(F_{C,n})^2},\nonumber
\end{align}
where in the last equality we used the product formula. 
Thus for any finite subset $\cM\subset \cM_{k}\backslash S_0$ we have
 \begin{align}\label{tolesszero}
\liminf_{n\to \infty}\sum_{v\in \cM_k\backslash S_0}N_v\cdot\frac{\log|\Res(F_{C, n})|_v}{\deg(F_{C,n})^2}&\leq \liminf_{n\to \infty} \sum_{v\in \cM}N_v\cdot\frac{\log|\Res(F_{C, n})|_v}{\deg(F_{C,n})^2}\nonumber\\
 &\leq -\sum_{v\in \cM}\log \capacity(M_{F,C,v})^{N_v}.
 \end{align}
Here for the last inequality we use the fact $1\leq \capacity(M_{F,C,v})\leq \displaystyle\liminf_{n\to\infty} |\Res F_{C, n}|_v^{-1/\deg(F_{C,n})^2}$. Combining \eqref{globalcapacitynonarch} and \eqref{tolesszero} we get that for any finite set $\cM\subset  \cM_k\backslash S_0$:
  $$\log \prod_{v\in \cM_k}\capacity(M_{F,C,v})^{N_v}\leq \sum_{v\in \cM_k\backslash S_0}\log \capacity(M_{F,C,v})^{N_v}-\sum_{v\in \cM}\log \capacity(M_{F,C,v})^{N_v}.$$
We may take an increasing sequence of finite sets $\cM_n\subset  \cM_k\backslash S_0$ such that $\displaystyle\cup_{n\geq 1}\cM_n=\cM_k\backslash S_0$, and apply the previous formula for $\cM_n$ in the place of $\cM$. Since the global capacity converges strongly and letting $n$ tend to $\infty$, inequality \eqref{negativeglobalcap} follows.
This finishes the proof of this proposition.
\end{proof}

\emph{Proof of theorem \ref{vchunlikely}.}
We combine Proposition \ref{equivalent heights} and Proposition \ref{samesmallpoints} to get 
\begin{align*}
\frac{\hat{h}_f(c)}{[k:\Q]}\cdot\log \prod_{v\in \cM_k} r_{in}(\mu_v)^{N_v}\leq \hat{h}_f(c)\cdot h(t)-\hhat_{f_{t}}(c(t))\leq \frac{\hat{h}_f(c)}{[k:\Q]}\cdot\log \prod_{v\in \cM_k} r_{out}(\mu_v)^{N_v},
\end{align*}
for all $t\in\overline{k}\setminus\Sing(f)$. Since the set $\Sing(f)$ is finite, we have that as $t\in\overline{k}$ varies
\begin{align*}
\hhat_{f_{t}}(c(t))=\hat{h}_f(c)h(t)+O(1),
\end{align*}
and the first part of our theorem follows.
The equidistribution statement in part \ref{thm3part1} now follows directly by combining Theorem \ref{arithmetic equidistribution} and Proposition \ref{samesmallpoints}.
\qed

\bibliographystyle{amsalpha}

\end{document}